\documentclass[final,3p]{elsarticle}  

\usepackage{bm}
\usepackage{amsmath}
\usepackage{amsfonts}
\usepackage{graphicx}
\usepackage{rotating}
\usepackage{multirow}
\usepackage[colorlinks,citecolor=black,filecolor=black,linkcolor=black,urlcolor=black]{hyperref}

\newcommand{\rot}[1]{\begin{sideways}#1\end{sideways}}

\begin{document}

\begin{frontmatter}

\title{A Sparse and High-Order Accurate Line-Based \\
Discontinuous Galerkin Method for Unstructured Meshes}

\author[UCBmath]{Per-Olof Persson\corref{cor1}}

\cortext[cor1]{Corresponding author. Tel.: +1-510-642-6947; Fax.: +1-510-642-8204.}
\ead{persson@berkeley.edu}
\address[UCBmath]{Department of Mathematics, University of California, Berkeley, Berkeley, CA 94720-3840, USA}

\begin{abstract}
We present a new line-based discontinuous Galerkin (DG) discretization scheme
for first- and second-order systems of partial differential equations.  The
scheme is based on fully unstructured meshes of quadrilateral or hexahedral
elements, and it is closely related to the standard nodal DG scheme as well as
several of its variants such as the collocation-based DG spectral element
method (DGSEM) or the spectral difference (SD) method. However, our motivation
is to maximize the sparsity of the Jacobian matrices, since this directly
translates into higher performance in particular for implicit solvers, while
maintaining many of the good properties of the DG scheme. To achieve this, our
scheme is based on applying one-dimensional DG solvers along each coordinate
direction in a reference element. This reduces the number of connectivities
drastically, since the scheme only connects each node to a line of nodes along
each direction, as opposed to the standard DG method which connects all nodes
inside the element and many nodes in the neighboring ones. The resulting
scheme is similar to a collocation scheme, but it uses fully consistent
integration along each 1-D coordinate direction which results in different
properties for nonlinear problems and curved elements. Also, the scheme uses
solution points along each element face, which further reduces the number of
connections with the neighboring elements. Second-order terms are handled by
an LDG-type approach, with an upwind/downwind flux function based on a switch
function at each element face. We demonstrate the accuracy of the method and
compare it to the standard nodal DG method for problems including Poisson's
equation, Euler's equations of gas dynamics, and both the steady-state and the
transient compressible Navier-Stokes equations. We also show how to integrate
the Navier-Stokes equations using implicit schemes and Newton-Krylov solvers,
without impairing the high sparsity of the matrices.
\end{abstract}

\end{frontmatter}

\section{Introduction}

In recent years it has become clear that the current computational methods for
scientific and engineering phenomena are inadequate for challenging
problems. These include problems with propagating waves, turbulent fluid flow,
nonlinear interactions, and multiple scales.  This has resulted in a
significant interest in so-called high-order accurate methods, which have the
potential to produce fundamentally more reliable solutions. A number of
numerical methods have been proposed, including multi-block finite difference
methods \cite{lele93compact,visbal02ho,nordstrom09multiblock}, high-order
finite volume methods \cite{Barth_kexact,gooch08hofvm}, stabilized finite
element methods \cite{hughes08stabilized}, discontinuous Galerkin (DG) methods
\cite{Reed_Hill,cockburn01rkdg,hesthaven08dgbook}, DG spectral element methods
(DGSEM) \cite{kopriva96}, spectral volume/difference methods
\cite{zj02spectralvolume,zj06spectraldifference,huynh07fluxreconstruction,vincent11fluxreconstruction},
and hybridized DG methods \cite{cockburn08hybrid,peraire09hybrid}. All methods
have advantages in particular situations, but for various reasons most general
purpose commercial-grade simulation tools still use traditional low-order
methods.

Much of the current research is devoted to the discontinuous Galerkin
method. This is partly because of its many attractive properties, including
the use of fully unstructured simplex meshes, the natural stabilization
mechanism based on approximate Riemann solvers, and the rigorous theoretical
foundations. It can certainly be discussed why the DG method is not used
routinely for real-world simulations, but one of the main reasons is clearly
its high computational cost, which is still at least a magnitude more than
low-order methods or high-order finite difference methods on similar
grids. For some problems, explicit time-stepping or matrix-free implicit
methods can be employed, but for many real-world problems and meshes full
Jacobian matrices are required for the solvers to be efficient. Here,
nodal-based Galerkin methods have a fundamental disadvantage in that they
connect all unknowns inside an element, as well as all neighboring face nodes,
even for first-order derivatives. This leads to a stencil size that scales
like $p^D$ for polynomial degrees $p$ in $D$ spatial dimensions. As a
contrast, a standard finite difference method only connects neighboring nodes
along the $D$ coordinate lines through the node. This gives a stencil size
proportional to $Dp$, which in three dimensions can be magnitudes smaller ever
for moderate values of $p$.

Several high-order schemes for unstructured meshes have been proposed with a
similar stencil-size reduction. In particular, the DG spectral element method
\cite{kopriva96,gassner10dgsem} is a collocation-based method on a staggered
grid which only uses information along each coordinate line for the
discretized equations.  Other closely related schemes have the same property,
such as the spectral difference method \cite{zj06spectraldifference}, the flux
reconstruction method
\cite{huynh07fluxreconstruction,vincent11fluxreconstruction}, and the DGM-FD
method \cite{hu11dgmfd}. For the special case of a linear one-dimensional
problem, many of these methods can be shown to be identical to the standard DG
method \cite{zj10unifying}, but in general they define different schemes with
varying properties.

In an attempt to further reduce the size of the Jacobians, and to ensure that
the scheme is identical to the standard DG method along each line of nodes, we
propose a new line-based DG scheme. Like the DGSEM, our Line-DG scheme is
derived by considering only the 1-D problems that arise along each coordinate
direction. We apply standard 1-D DG formulations for each of these
sub-problems, and all integrals are computed fully consistently (with
sufficient accuracy), which means in particular that the definition of the
scheme makes no statement about flux points. We note that this can be done
without introducing additional connectivities, since all nodes in the local
1-D problem are already connected by the shape functions. In addition, our
scheme uses solution points along each element face, which further reduces the
number of connectivities with the neighboring elements.

For the second-order terms in the Navier-Stokes equations, we use an LDG-type
approach \cite{cockburn98ldg} with upwind/downwind fluxes based on consistent
switches along all globally connected lines of elements. Special care is
required to preserve the sparsity of the resulting matrices, and we propose a
simple but efficient Newton-Krylov solver which splits the matrix product in
order to avoid introducing additional matrix entries. Many options for
preconditioning are possible, and in this work we use a block-Jacobi method
with sparse blocks.

We first describe the method for first-order systems in Section 2 and mention
some practical implementation issues, including a study of the structure of
the Jacobian matrices. In Section 3 we extend the scheme to second-order
systems using the LDG-type scheme, and in Section 4 we discuss the implicit
temporal discretization, some approaches for maintaining the high sparsity of
the discretization, and the Newton-Krylov solver. Finally, in Section 5 we
show numerical results and convergence for Poisson's equation, an inviscid
Euler vortex, and flow over a cylinder. We also compare the method to the
standard nodal DG method, and we conclude that the differences are overall
very small. For the Navier-Stokes equations, we show convergence of drag and
lift forces for steady-state laminar flow around an airfoil, and we
demonstrate our implicit time-integrators on a transient LES-type flow
problem.

\section{Line-based discontinuous Galerkin discretization}

\subsection{First-order equations}

\begin{figure}
  \includegraphics[width=\textwidth]{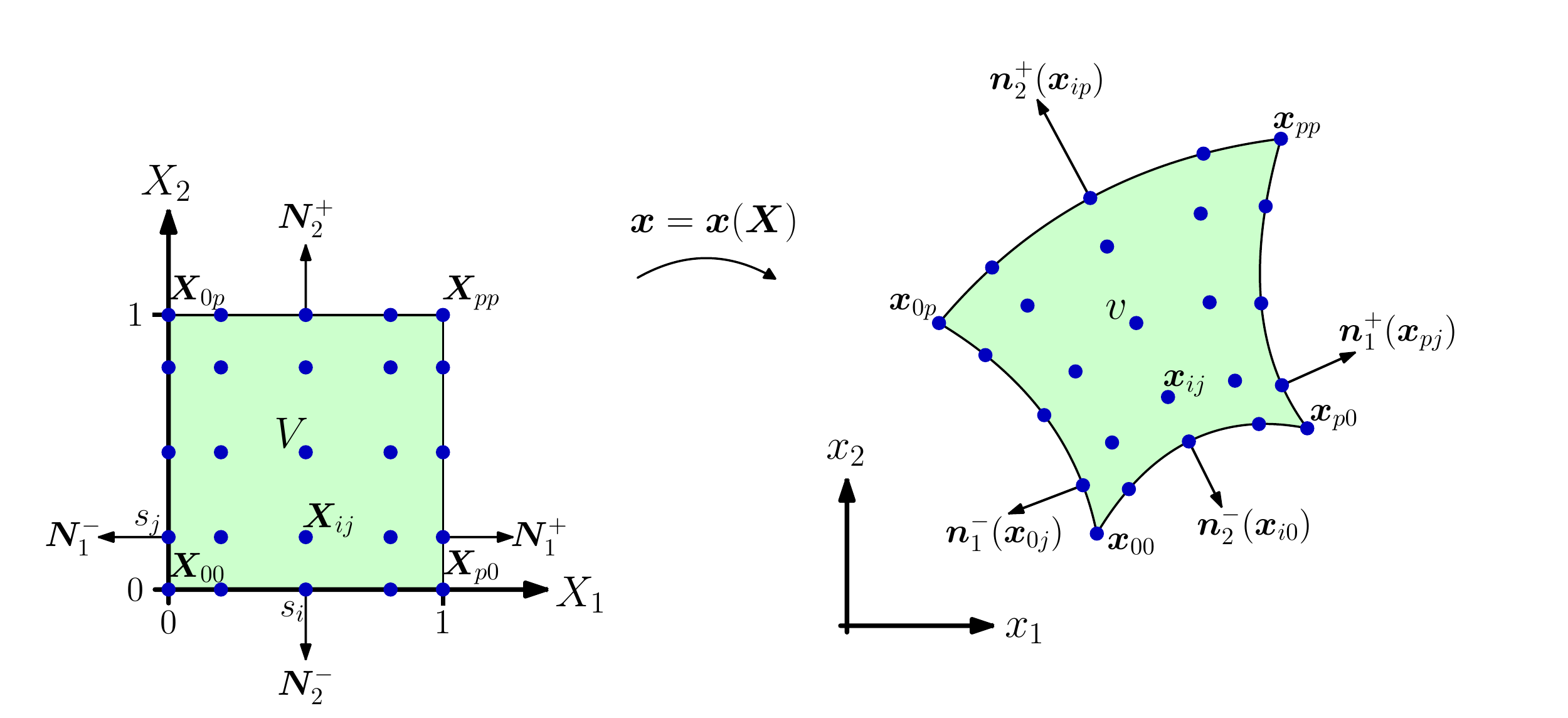}
  \caption{A two-dimensional illustration of the mapping from a reference
           element $V$ to the actual curved element $v$, for the case $p=4$.}
  \label{plt}
\end{figure}

Consider a system of $m$ first-order conservation laws with source terms,
\begin{align}
\frac{\partial \bm{u}}{\partial t} + \nabla \cdot \bm{F}(\bm{u}) = \bm{S}(\bm{u}), \label{conslaw}
\end{align}
in a three-dimensional domain $\Omega$, with solution $\bm{u}$, flux function
$\bm{F}(\bm{u})$, source function $\bm{S}(\bm{u})$, and appropriate boundary
conditions on $\partial \Omega$. We will use a discretization of $\Omega$ into
non-overlapping, conforming, curved hexahedral elements. Within each element
we introduce a Cartesian grid of $(p+1)^3$ node points, where $p\ge 1$, by
defining a smooth one-to-one mapping given by a diffeomorphism
$\bm{x}=\bm{x}(\bm{X})$ between the reference unit cube $V=[0,1]^3$ and the
element $v$, and setting $\bm{x}_{ijk}=\bm{x}(\bm{X}_{ijk})$, where
$\bm{X}_{ijk}=(s_i,s_j,s_k)$ for $0\le i,j,k\le p$, and $\{s_i\}$ is an
increasing sequence of $p+1$ node positions $s_i\in[0,1]$ with $s_0=0$ and
$s_p=1$ (see figure~\ref{plt}).

To obtain our numerical scheme for approximating (\ref{conslaw}), we
consider a single element $v$ and its mapping $\bm{x}=\bm{x}(\bm{X})$,
and follow standard procedure to change independent variables from
$\bm{x}$ to $\bm{X}$. This transforms (\ref{conslaw}) into
\begin{align}
J\frac{\partial \bm{u}}{\partial t} + \nabla_{\bm{X}} \cdot \widetilde{\bm{F}}(\bm{u}) = J \bm{S}(\bm{u}), \label{conslawref}
\end{align}
in the reference domain $V$. Here we have defined the mapping Jacobian
$J=\det(\bm{G})$ and the contravariant fluxes $\widetilde{\bm{F}} =
(\widetilde{\bm{f}}_1,\widetilde{\bm{f}}_2,\widetilde{\bm{f}}_3 ) = J
\bm{G}^{-1} \bm{F}$, with the mapping deformation gradient
$\bm{G}=\nabla_{\bm{X}} \bm{x}$.

A standard nodal discontinuous Galerkin method would now consider the
multivariate polynomial $\bm{u}(\bm{X})$ that interpolates the grid function,
$\bm{u}_{ijk} = \bm{u}(\bm{X}_{ijk})$, and define a numerical scheme for the
spatial derivatives of (\ref{conslaw}) by a Galerkin procedure in $V$. Our
approach differs in that it considers each of the three spatial derivatives in
(\ref{conslawref}) separately and approximates them numerically using
one-dimensional discontinuous Galerkin formulations along each of the $3(p+1)^2$
curves defined by straight lines in the reference domain $V$, through the
three sets of nodes along each space dimension.

More specifically, the $(p+1)^2$ curves along the first space dimension are
$\bm{x}_{jk}(\xi)=\bm{x}(\xi,X_j,X_k)$ for $0\le j,k \le p$. On these we
define the polynomial $\bm{u}_{jk}(\xi)\in \mathcal{P}_p([0,1])^m$ that
interpolates $\bm{u}_{ijk}$, $i=0,\ldots,p$, and we define a numerical
approximation $\bm{r}_{jk}(X_1)$ to $\partial \widetilde{\bm{f}_1}/\partial
X_1$ by a one-dimensional Galerkin procedure: Find
$\bm{r}_{jk}(\xi)\in\mathcal{P}_p([0,1])^m$ such that
\begin{align}
&\int_0^1 \bm{r}_{jk}(\xi) \cdot \bm{v}(\xi)\, d\xi = \int_0^1 \frac{d\widetilde{\bm{f}_1}(\bm{u}_{jk}(\xi))}{d\xi}\cdot \bm{v}(\xi)\,d\xi \nonumber \\
&\qquad = \widehat{\widetilde{\bm{f}_1}}(\bm{u}_{jk}^+(1), \bm{u}_{jk}(1)) \cdot \bm{v}(1) - 
   \widehat{\widetilde{\bm{f}_1}}(\bm{u}_{jk}(0), \bm{u}_{jk}^-(0)) \cdot \bm{v}(0) -
   \int_0^1 \widetilde{\bm{f}_1}(\bm{u}_{jk}(\xi)) \cdot \frac{d\bm{v}}{d\xi}\,d\xi, \label{galerkinscheme}
\end{align}
for all test functions $\bm{v}(\xi)\in\mathcal{P}_p([0,1])^m$. Here,
$\bm{u}_{jk}^+(1)$ is the numerical solution at
$\bm{x}_{jk}(1^+)=\bm{x}(1^+,X_j,X_k)$, and similarly
$\bm{u}_{jk}^-(0)$ at $\bm{x}_{jk}(0^-)=\bm{x}(0^-,X_j,X_k)$. These
will be given either by nodes in the neighboring elements or
implicitly through the boundary conditions. Furthermore,
$\widehat{\widetilde{\bm{f}_1}}(\bm{u}_R,\bm{u}_L)$ is a numerical
flux function for $\widetilde{\bm{f}_1}$, but we note that with the
reference normal direction $\bm{N}_1^+=(1,0,0)$, the contravariant
flux can be written
\begin{align}
\widetilde{\bm{f}}_1 = \widetilde{\bm{F}} \cdot \bm{N}_1^+
= (J\bm{G}^{-1}\bm{F})\cdot \bm{N}_1^+ = \bm{F}\cdot (J\bm{G}^{-T}\bm{N}_1^+)
= \bm{F}\cdot\bm{n}_1^+ \label{contraflux}
\end{align}
with the (non-normalized) normal vector $\bm{n}_1^+=J\bm{G}^{-T}\bm{N}_1^+$ at
the boundary point $\bm{x}_{jk}(1)$. Our numerical flux then becomes
\begin{align}
\widehat{\widetilde{\bm{f}_1}}(\bm{u}_R,\bm{u}_L) =
\widehat{\widetilde{\bm{F}}\cdot\bm{N}_1^+}(\bm{u}_R,\bm{u}_L) =
\widehat{\bm{F}\cdot\bm{n}_1^+} (\bm{u}_R,\bm{u}_L), \label{numflux1}
\end{align}
where $\widehat{\bm{F}\cdot\bm{n}}(\bm{u}^+,\bm{u}^-)$ is a standard numerical
flux function used in finite volume and discontinuous Galerkin schemes, with
normal direction $\bm{n}$ and traces $\bm{u}^\pm$ in the positive/negative
normal direction. This allows us to use existing flux functions and
approximate Riemann solvers without modification.

Similarly, for the second numerical flux we move the negative sign to
the normal direction and define $\bm{N}_1^-=(-1,0,0)$ and
$\bm{n}_1^-=J\bm{G}^{-T}\bm{N}_1^-$, which is again an outward normal vector.
We can then write:
\begin{align}
-\widehat{\widetilde{\bm{f}_1}}(\bm{u}_R,\bm{u}_L) =
\widehat{\widetilde{\bm{F}}\cdot\bm{N}_1^-}(\bm{u}_L,\bm{u}_R) =
\widehat{\bm{F}\cdot\bm{n}_1^-} (\bm{u}_L,\bm{u}_R), \label{numflux2}
\end{align}
where we also have swapped the order of the arguments to the flux function
$\widehat{\bm{F}\cdot\bm{n}}$ to be consistent with the negative normal
direction.  Our Galerkin scheme (\ref{galerkinscheme}) then gets the final
form: Find $\bm{r}_{jk}(\xi)\in\mathcal{P}_p([0,1])^m$ such that
\begin{align}
&\int_0^1 \bm{r}_{jk}(\xi) \cdot \bm{v}(\xi)\, d\xi =
\widehat{\bm{F}\cdot\bm{n}_1^+} (\bm{u}_{jk}^+(1), \bm{u}_{jk}(1)) \cdot \bm{v}(1) +
\widehat{\bm{F}\cdot\bm{n}_1^-} (\bm{u}_{jk}^-(0), \bm{u}_{jk}(0)) \cdot \bm{v}(0) -
   \int_0^1 \widetilde{\bm{f}_1}(\bm{u}_{jk}(\xi)) \cdot \frac{d\bm{v}}{d\xi}\,d\xi, \label{galerkinscheme2}
\end{align}
for all $\bm{v}(\xi)\in\mathcal{P}_p([0,1])^m$.

We use a standard finite element procedure to solve (\ref{galerkinscheme2})
for $\bm{r}_{jk}(\xi)$.  Introduce the nodal Lagrange basis functions
$\phi_i\in\mathcal{P}_p([0,1])$ such that $\phi_i(s_j)=\delta_{ij}$, for
$i,j=0,\ldots,p$, and set
\begin{align}
\bm{u}_{jk}(\xi) &= \sum_{i=0}^p \bm{u}_{ijk} \phi_i(\xi), \\
\bm{r}_{jk}(\xi) &= \sum_{i=0}^p \bm{r}_{ijk} \phi_i(\xi), \label{reqn}
\end{align}
To find the $m(p+1)$ coefficients along the curve $\bm{x}_{jk}(\xi)$, we set
$\bm{v}(\xi)=\bm{e}_\ell \phi_i(\xi)$, for each $i=0,\ldots,p$ and
$\ell=1,\ldots,m$, where $(\bm{e}_\ell)_n = \delta_{\ell n}$ for
$n=1,\ldots,m$. Our Galerkin scheme (\ref{galerkinscheme2}) then gets the
discrete form $\bm{M} \bm{r}_{jk} = \bm{b}$, and we find the coefficients
$\bm{r}_{jk}$ by solving $m$ linear systems with the $(p+1)$-by-$(p+1)$ mass matrix
$\bm{M}$. Repeating the procedure for each $j,k=0,\ldots,p$ we obtain all
coefficients $\bm{r}_{ijk}=\bm{r}^{(1)}_{ijk}$, which is the grid function for
our numerical approximation of $\partial \widetilde{\bm{F}}_1/\partial X_1$ at
each grid point $\bm{x}_{ijk}$.

In an analogous way, we calculate coefficients $\bm{r}_{ijk}^{(2)}$
and $\bm{r}_{ijk}^{(3)}$ that approximate $\partial
\widetilde{\bm{F}}_2/\partial X_2$ and $\partial
\widetilde{\bm{F}}_3/\partial X_3$, respectively, at the grid points.
The curves considered are now $\bm{x}_{ik}=\bm{x}(X_i,\xi,X_k)$ and
$\bm{x}_{ij}=\bm{x}(X_i,X_j,\xi)$, and with the reference normals
$\bm{N}_2^{\pm}=(0,\pm 1,0)$ and $\bm{N}_3^{\pm}=(0,0,\pm 1)$ the
contravariant fluxes $\widetilde{\bm{f}}_2$ and $\widetilde{\bm{f}}_3$
can again be written as $\bm{F}\cdot\bm{n}$ where
$\bm{n}=J\bm{G}^{-T}\bm{N}$ is a non-normalized normal vector to the
element at the boundary points. The solution procedure involves the
same mass matrix $\bm{M}$ and is identical to before.

Using the calculated numerical approximations to each partial derivative
in (\ref{conslawref}), we obtain our final semi-discrete formulation:
\begin{align}
\frac{d\bm{u}_{ijk}}{dt} + \frac{1}{J_{ijk}}\sum_{n=1}^3 \bm{r}_{ijk}^{(n)} = \bm{S}(\bm{u}_{ijk}), \label{semidisc}
\end{align}
where $J_{ijk}=J(\bm{x}_{ijk})$. 

\subsection{Implementation details}

For the mapping $\bm{x}(\bm{X})$ it is natural to use an iso-parametric
approach. The node positions $\bm{x}_{ijk}$ are given by some curved mesh
generation procedure \cite{persson09curved}, and we define
\begin{align}
\bm{x}(\bm{X})=\sum_{i,j,k=0}^p \bm{x}_{ijk} \phi_i(X_1)\phi_j(X_2)\phi_k(X_3),
\end{align}
which clearly satisfies our interpolation requirement
\begin{align}
\bm{x}(\bm{X}_{ijk}) &= \sum_{i',j',k'=0}^p \bm{x}_{i'j'k'} \phi_{i'}(s_i)\phi_{j'}(s_j)\phi_{k'}(s_k) 
 =\sum_{i',j',k'=0}^p \bm{x}_{i'j'k'} \delta_{ii'}\delta_{jj'}\delta_{kk'} = \bm{x}_{ijk}.
\end{align}
This allows us to easily compute $\bm{G}(\bm{X})$ at any point $\bm{X}$, which
will involve the derivatives $\phi_i'(\xi)$ of the shape functions. To
evaluate the one-dimensional integrals in (\ref{galerkinscheme2}), we use
Gauss-Legendre integration of sufficiently high degree. For all our problems,
a precision of $3p$ appears to be enough, so we use integration rules with
$\lceil (3p+1)/2 \rceil$ integration points.

The computation of the discretization (\ref{galerkinscheme2}) is remarkably
simple compared to a nodal DG scheme, primarily because (a) The integrals are
only one-dimensional, and (b) The numerical fluxes are only evaluated
point-wise.  We note that the left-hand side of (\ref{galerkinscheme2}), which
contributes to the mass matrix $\bm{M}$, is constant regardless of solution
component, line, and element (even if the actual mapped element is
curved). Therefore it can be pre-computed and pre-factorized using a standard
Cholesky method. Furthermore, many lines and components can be processed
simultaneously, which might further increase the performance through the use
of BLAS3-type cache-optimized linear algebra libraries.

For the integral in the right-hand side of (\ref{galerkinscheme2}), the term
$d\bm{v}/d\xi$ is again constant for all components, lines, and elements, so
its discretization at the Gauss integration points can be pre-computed and
combined with the inverted mass matrix and the Gauss integration weights
$\bm{w}$. For non-linear problems, the only part that requires re-evaluation
at each Gauss integration point is $\widetilde{\bm{f}}_1(\bm{u}_{jk}(\xi))$,
although the deformation gradient $\bm{G}$ can be pre-computed if necessary.

For the numerical fluxes, we pointed out above that (\ref{numflux1}) and
(\ref{numflux2}) have exactly the same form as standard numerical flux
functions. We pre-compute the outward normals $\bm{n}_i^+$ and $\bm{n}_i^-$,
for $i=1,2,3$, at all boundary nodes. Note that our scheme only computes
point-wise numerical fluxes, unlike the nodal DG method which involves
integrals of the numerical fluxes. Sometimes existing numerical flux functions
require the normal vector to be of unit length, in this case we normalize
$\bar{\bm{n}}=\bm{n}/|\bm{n}|$ and use the fact that
\begin{align}
\widehat{\bm{F}\cdot\bm{n}} = |\bm{n}| \widehat{\bm{F}\cdot\bar{\bm{n}}}.
\end{align}

Finally, the multipliers $J_{ijk}$ are defined at the node points (not the
Gauss integration points), and can also be pre-computed.

\subsection{Stencil size and sparsity pattern}

To illustrate the drastic reduction of the number of entries in the Jacobian
matrices for the Line-DG method, consider the $(p+1)^3$ nodes in an (interior)
element and its six neighboring elements. For a first-order operator, we note
that a standard nodal DG formulation will in general produce full block
matrices, that is, each degree of freedom will depend on all the other ones
within the element. In addition, the face integrals will connect all nodes on
an element face to all neighboring element face nodes.  This gives
$6(p+1)^2(p+1)^2=6(p+1)^4$ additional connections per element, or in average
$6(p+1)^4/(p+1)^3=6(p+1)$ connections per degree of freedom. In total, the
average number of connections is $(p+1)^3+6(p+1)$, which illustrates why
matrix-based DG methods are considered memory intensive and expensive even at
modest values of $p$.

\begin{figure}
  \begin{center}
    \includegraphics[width=.75\textwidth]{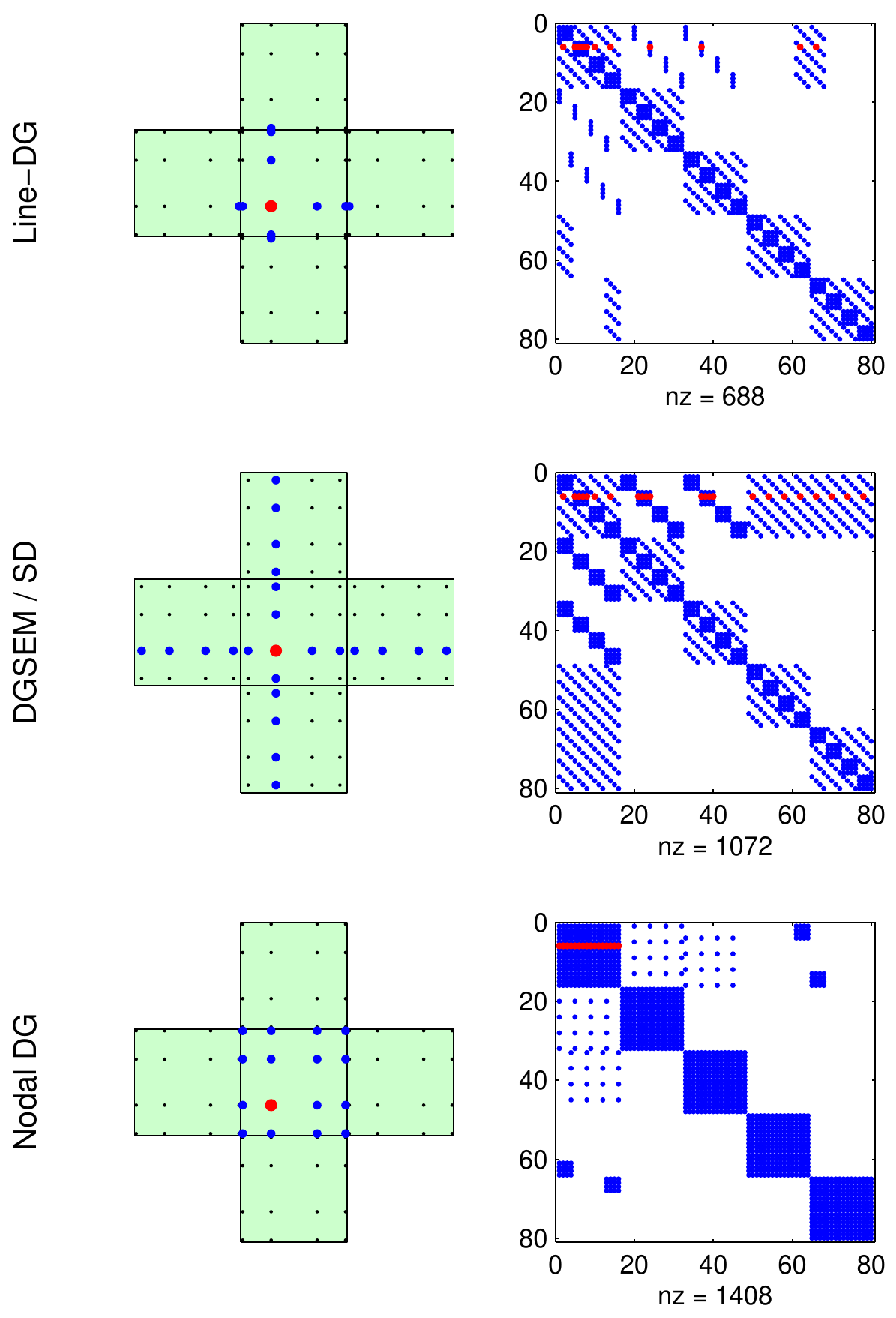}
  \end{center}
  \caption{The connectivities (blue circles) to a single node (red circle)
    for the Line-DG method, the DGSEM/SD method, and the nodal DG method
    (2-D quadrilateral elements, a first-order operator).}
  \label{fig1}
\end{figure}

As a contrast, in our line-based method each node will only connect to other
nodes within the same lines, and to only one node in each neighboring element,
for a total of $(3p+1)+6=3p+7$ connectivities. This is similar to that of the
DGSEM/SD methods, although with Gauss-Legendre solution points these schemes
also connect entire lines of nodes in the neighboring element, giving a total
of $(3p+1)+6(p+1)=9p+7$ connectivities.  These numbers are tabulated for a
range of degrees $p$ in three dimensions in table~\ref{costtab}. The sparsity
patterns are illustrated in figure~\ref{fig1} for two-dimensional
quadrilateral elements, for all three methods. The connectivities are shown
both by a nodal plot, with bold nodes corresponding to the dependencies of the
single red node, and by sparsity plots of the Jacobian matrices.

We note that in three dimensions, already for $p=3$ the Line-DG method is 5.5 times
sparser than nodal DG, and for $p=10$ it is almost 40 times sparser. This
reduction in stencil size translates into lower assembly times, but more
importantly, for matrix-based solvers it means drastically lower storage
requirements and faster matrix-vector products for iterative implicit solvers.

\begin{table}
\begin{center}
\begin{tabular}{ll|rrrrrrrrrr}
\hline
& Polynomial order $p$ & $1$ & $2$ & $3$ & $4$ & $5$ & $6$ & $7$ & $8$ & $9$ & $10$ \\
\hline
\multirow{3}{*}{\rot{2-D}} & Line-DG connectivities & $7$ & $9$  & $11$  & $13$ & $15$  & $17$ & $19$ &  $21$ & $23$ & $25$ \\
 & DGSEM/SD connectivities & $11$ & $17$  & $23$  & $29$ & $35$  & $41$ & $47$ &  $53$ & $59$ & $65$ \\
 & Nodal DG connectivities & $8$ & $13$  & $20$  & $29$ & $40$  & $53$ & $68$ &  $85$ & $104$ & $125$ \\ \hline
\multirow{3}{*}{\rot{3-D}} & Line-DG connectivities & $10$ & $13$  & $16$  & $19$ & $22$  & $25$ & $28$ &  $31$ & $34$ & $37$ \\
 & DGSEM/SD connectivities & $16$ & $25$  & $34$  & $43$ & $52$  & $61$ & $70$ &  $79$ & $88$ & $97$ \\
 & Nodal DG connectivities & $20$ & $45$  & $88$  & $155$ & $252$  & $385$ & $560$ &  $783$ & $1060$ & $1397$ \\ \hline
\end{tabular}
\end{center}
\caption{The number of connectivities per node for a first-order operator and 2-D quadrilateral / 3-D hexahedral elements with the Line-DG, the DGSEM/SD, and the nodal DG methods.}
\label{costtab}
\end{table}

\section{Second-order equations}

We now consider the discretization of equations with second-order
derivatives, in the form of a system of conservation laws
\begin{align}
\frac{\partial \bm{u}}{\partial t} +
\nabla \cdot \bm{F}(\bm{u}, \nabla{\bm{u}}) = \bm{S}(\bm{u},\nabla\bm{u}).
\label{system1}
\end{align}
We first use a standard technique in many finite difference and
discontinuous Galerkin methods, and introduce the auxiliary variables
$\bm{q}$ and rewrite as a split system
\begin{align}
\frac{\partial \bm{u}}{\partial t} +
\nabla \cdot \bm{F}(\bm{u}, \bm{q}) &= \bm{S}(\bm{u},\bm{q}), \label{spliteq1} \\
  \nabla \bm{u} &= \bm{q} . \label{spliteq2}
\end{align}
This essentially has the form of our first-order system (\ref{conslaw}), and
we can apply Line-DG to each solution component as described above. More
specifically, the change of variables from $\bm{x}$ to $\bm{X}$ transforms
(\ref{spliteq1}), (\ref{spliteq2}) into
\begin{align}
J \frac{\partial\bm{u}}{\partial t} + \nabla_{\bm{X}} \cdot
\widetilde{\bm{F}}(\bm{u},\bm{q}) &= J \bm{S}(\bm{u},\bm{q}), \label{splitdisc1} \\
\nabla_{\bm{X}} \cdot \widetilde{\bm{u}}(\bm{u}) &= J \bm{q}, \label{splitdisc2}
\end{align}
where $\widetilde{\bm{u}} =
(\widetilde{\bm{u}}_1,\widetilde{\bm{u}}_2,\widetilde{\bm{u}}_3) = \bm{u} \otimes J \bm{G}^{-1}$. We discretize (\ref{splitdisc1}) as described before
for the first-order case, treating $\bm{q}$ as additional solution
components. For (\ref{splitdisc2}), we use a completely analogous
procedure. We introduce the grid function
$\bm{q}_{ijk}=\bm{q}(\bm{X}_{ijk})$, and along each curve
$\bm{X}_{jk}(\xi)$ we define the polynomial
$\bm{q}_{jk}(\xi)\in\mathcal{P}_p([0,1])^{m\times 3}$ that
interpolates $\bm{q}_{ijk}$, $i=0,\ldots,p$. We find the numerical
approximation $\bm{d}_{jk}(X_1)$ to $\partial
\widetilde{\bm{u}}_1/\partial X_1$ by the Galerkin formulation: Find
$\bm{d}_{jk}(\xi)\in\mathcal{P}_p([0,1])^{m\times 3}$ such that
\begin{align}
&\int_0^1 \bm{d}_{jk}(\xi) : \bm{\tau}(\xi)\,d\xi = 
\int_0^1 \frac{d \widetilde{\bm{u}}_1}{d\xi} : \bm{\tau}(\xi)\,d\xi \nonumber \\
&\qquad\qquad=\widehat{\widetilde{\bm{u}}}_1(\bm{u}_{jk}^+(1),\bm{q}_{jk}^+(1),\bm{u}_{jk}(1),\bm{q}_{jk}(1))
   : \bm{\tau}(1) - 
 \widehat{\widetilde{\bm{u}}}_1(\bm{u}_{jk}(0),\bm{q}_{jk}(0),\bm{u}_{jk}^-(0),\bm{q}_{jk}^-(0))
   : \bm{\tau}(0) \nonumber \\
& \qquad \qquad\ \ \ \ %
- \int_0^1 \widetilde{\bm{u}}_1(\bm{u}_{jk}(\xi)) : \frac{d\bm{\tau}}{d\xi}\,d\xi
\end{align}
for all test functions $\bm{\tau}(\xi)\in\mathcal{P}_p([0,1])^{m\times 3}$.
Note that we allow for the numerical flux $\widehat{\widetilde{\bm{u}}}_1$ to
depend on both $\bm{u}$ and $\bm{q}$ on each side of the face, even though
the actual flux $\widetilde{\bm{u}}_1$ is only a function of $\bm{u}$. Again,
the numerical contravariant fluxes can be written in terms of the actual
fluxes and the actual normal vector:
\begin{align}
\widehat{\widetilde{\bm{u}_1}} =
\widehat{\widetilde{\bm{u}}\cdot\bm{N}_1^+} =
\widehat{\bm{u}\otimes \bm{n}_1^+} = \widehat{\bm{u}}\otimes\bm{n}_1^+, \label{numfluxu1}
\end{align}
and similarly in the negative direction and along the other coordinate
directions.  It remains only to define the numerical fluxes
$\widehat{\bm{F}\cdot\bm{n}}=\widehat{\bm{F}}\cdot\bm{n}$ and
$\widehat{\bm{u}}$. We could in principle consider any scheme that can be
written in this form \cite{arnold02unified}, such as the interior penalty
method, the BR2 method, the LDG method \cite{cockburn98ldg}, and the CDG
method \cite{peraire08cdg}. Here we use a scheme based on the LDG method,
because it has a simple upwind/downwind character, it does not evaluate
derivatives of grid functions at the boundaries, and it appears well-suited
for our Line-DG discretization. Furthermore, since our implicit solvers avoid
the elimination of $\bm{q}$, the scheme has a compact connectivity (only
connects neighboring elements).

First, we separate the fluxes $\bm{F}$ into an inviscid and a viscous
part:
\begin{align}
\bm{F}(\bm{u},\nabla \bm{u}) = \bm{F}^\mathrm{inv}(\bm{u}) + \bm{F}^\mathrm{vis}(\bm{u},\nabla \bm{u}).
\end{align}
This decomposition is clearly not unique, but it is understood that
for many problems there is a natural separation into a
convection-dominated inviscid component and a diffusion-dominated
viscous component. This allows us to use standard approximate Riemann
solvers for $\widehat{\bm{F}}^\mathrm{inv}$ as before, and we will now
consider only the treatment of the viscous fluxes
$\widehat{\bm{F}}^\mathrm{vis}$. For shorter notation, we assume below
that $\bm{F}=\bm{F}^\mathrm{vis}$ and that $\bm{n}$ is a unit vector.

We will define the fluxes in terms of a so-called switch function, which
simply assigns a sign to each internal element face. For one-dimensional
problems, the natural switch function is to set all these signs equal (either
positive or negative), and we will mimic this for our Line-DG method by
identifying globally connected lines in our hexahedral meshes.

\begin{figure}
  \begin{center}
    \includegraphics[width=.5\textwidth]{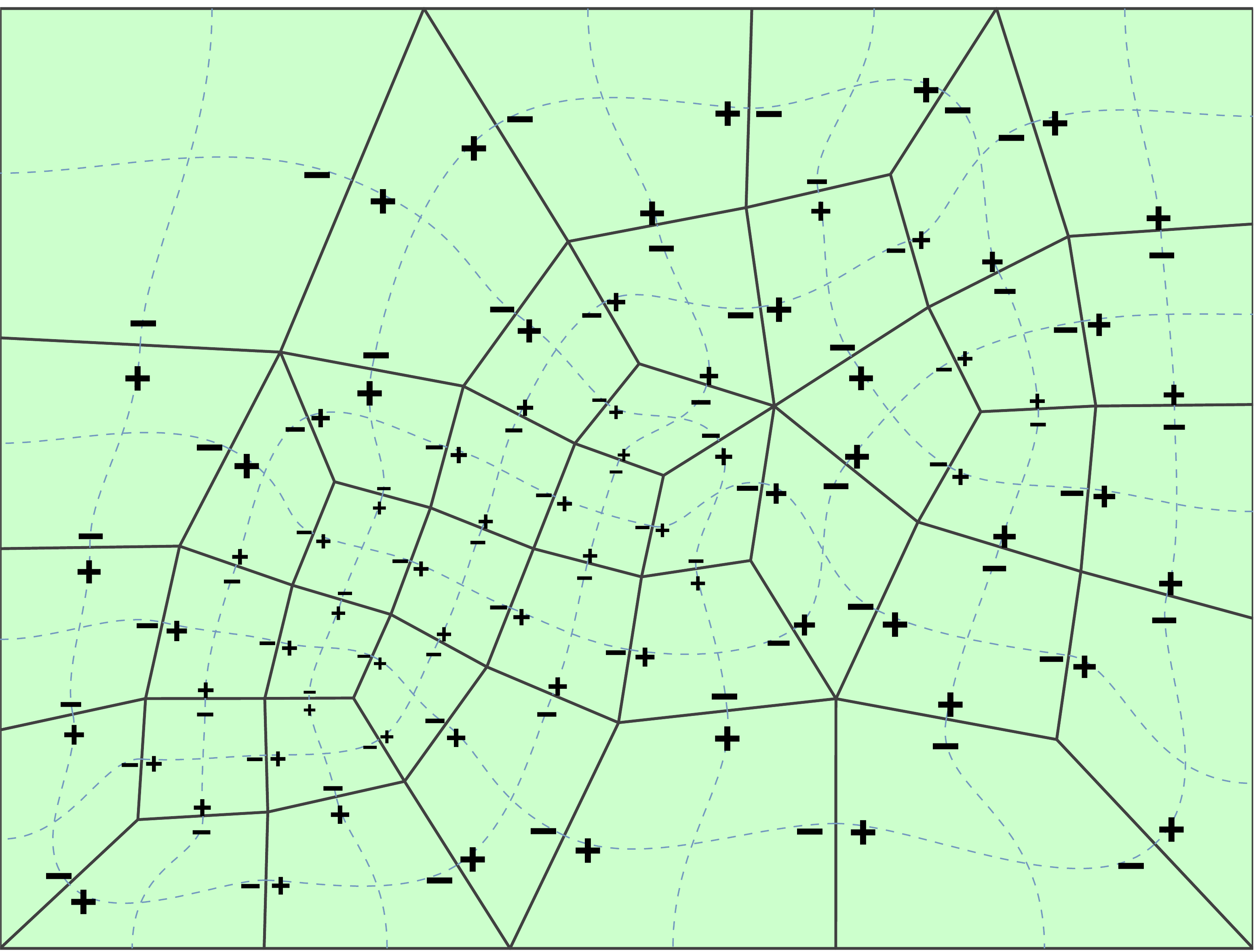}
  \end{center}
  \caption{A sample quadrilateral mesh and switch function $S^\pm_i$ for
    $i=1,2$ in each element. Note that the switches have consistent directions
    along each one-dimensional global curve (dashed blue lines).}
  \label{switchplt}
\end{figure}

In our notation, instead of assigning switches to each face, we introduce
$S^{\pm}_i \in \{-1,1\}$ for the switches at local coordinate $\xi=1$ and
$\xi=0$ along direction $i=1,2,3$. There is some redundancy here, since we
require that $S^+_i = -S^-_i$, and also that the switch function for a shared
face between two neighboring elements have opposite signs. See
figure~\ref{switchplt} for an example quadrilateral mesh and switch function.
This was generated by a straight-forward algorithm, where an arbitrary element
face is chosen and assigned an arbitrary sign, which then defines the
alternating pattern along a sequence of elements in both directions. This
procedure is repeated until all faces have been processed.

With the switch function defined, we can formulate the LDG numerical fluxes
for the second-order terms:
\begin{align}
\widehat{\bm{F}}(\bm{u},\bm{q},\bm{n}) &= 
\{\!\{ \bm{F}(\bm{u},\bm{q}) \}\!\} +
C_{11} [\![ \bm{u} \otimes \bm{n} ]\!] + \bm{C}_{12} \otimes [\![ \bm{F}(\bm{u},\bm{q})\cdot\bm{n} ]\!], \label{ldgflux1} \\
\widehat{\bm{u}}(\bm{u},\bm{q},\bm{n}) &=
\{\!\{ \bm{u} \}\!\} -\bm{C}_{12}\cdot [\![ \bm{u}\otimes\bm{n} ]\!] + C_{22}[\![\bm{F}(\bm{u},\bm{q})\cdot\bm{n}]\!]. \label{ldgflux2}
\end{align}
for a solution $\bm{u},\bm{q}$ and a face normal vector $\bm{n}$.  Here,
$\{\!\{ \cdot \}\!\}$ denotes the mean value and $[\![ \cdot ]\!]$ denotes the
jump over a face:
\begin{align}
\{\!\{ \bm{v} \}\!\} \equiv \frac12 (\bm{v}^++\bm{v}^-),\qquad
[\![ \bm{v} \odot \bm{n} ]\!] \equiv \bm{v}^+\odot\bm{n}-\bm{v}^-\odot\bm{n}
\end{align}
where $\bm{v}^+$ is the quantity $\bm{v}$ on the positive side of the face
(according to the normal $\bm{n}$), $\bm{v}^-$ is $\bm{v}$ on the negative
side, and $\odot$ is any multiplication operator. The coefficients
$C_{11},\bm{C}_{12},C_{22}$ give the scheme different properties, and we note
in particular that:
\begin{itemize}
\item If $C_{22}=0$, the fluxes (\ref{ldgflux2}) do not depend on $\bm{q}$, which means
      the discretized equation (\ref{splitdisc2}) immediately give $\bm{q}$ within each
      element (no coupling to equation (\ref{splitdisc1})).
\item For the particular choice $\bm{C}_{12} = \bm{n} S^\pm_i /2$, where $S^\pm_i$ is
      the switch for the considered face and direction, the fluxes can
      be written in the form:
\begin{align}
\widehat{\bm{F}}(\bm{u}_R,\bm{q}_R,\bm{u}_L,\bm{q}_L,\bm{n}) &= 
C_{11} [\![ \bm{u} \otimes \bm{n} ]\!] + 
\begin{cases}
\bm{F}(\bm{u}_R,\bm{q}_R) & \text{if } S^\pm_i=+1 \\
\bm{F}(\bm{u}_L,\bm{q}_L) & \text{if } S^\pm_i=-1 \\
\end{cases} \\
\widehat{\bm{u}}(\bm{u}_R,\bm{q}_R,\bm{u}_L,\bm{q}_L,\bm{n}) &= C_{22}[\![\bm{F}(\bm{u},\bm{q})\cdot\bm{n}]\!] +
\begin{cases}
\bm{u}_L & \text{if } S^\pm_i=+1 \\
\bm{u}_R & \text{if } S^\pm_i=-1 \\
\end{cases}
\end{align}
where it is clear how the method is upwinding/downwinding the two numerical
fluxes, depending on the switch $S^\pm_i$. The constants $C_{11}$ and $C_{22}$ are
additional stabilization parameters, which can be seen as penalties on the
jumps in the solution and in the normal fluxes, respectively. In many of our
problems we set both of these coefficients to zero (the so-called minimal
dissipation LDG method \cite{cockburn07mdldg}). This makes the scheme
particularly simple, and also further reduces the number of connectivities in
the Jacobian matrices.
\end{itemize}
At the boundaries we impose conditions by appropriate choices of numerical
fluxes. For example, at a Dirichlet-type boundary with a prescribed solution
$\bm{u}=\bm{g}_D$, we set:
\begin{align}
\widehat{\bm{F}}(\bm{u},\bm{q},\bm{n}) &= 
\bm{F}(\bm{u},\bm{q}) +
C_{11} ( \bm{u} - \bm{g}_D ) \otimes \bm{n} \label{dirichlet1} \\
\widehat{\bm{u}}(\bm{u},\bm{q},\bm{n}) &= \bm{g}_D, \label{dirichlet2}
\end{align}
where $C_{11}$ in (\ref{dirichlet1}) must be positive, even though we often
choose $C_{11}=0$ for the interior fluxes. At a Neumann-type boundary with
prescribed normal fluxes $\bm{F}\cdot\bm{n}=\bm{g}_N$, we set:
\begin{align}
\widehat{\bm{F}}(\bm{u},\bm{q},\bm{n}) &= \bm{g}_N \otimes \bm{n} \label{neumann1} \\
\widehat{\bm{u}}(\bm{u},\bm{q},\bm{n}) &= \bm{u} -C_{22} (\bm{F}(\bm{u},\bm{q})\cdot\bm{n}-\bm{g}_N). \label{neumann2}
\end{align}
For mixed conditions we apply combinations of these fluxes for the different
components of $\bm{u}$ and $\bm{q}$. 

With the fluxes defined, we can calculate
$\bm{d}_{ijk}=\bm{d}^{(1)}_{ijk}$ for all $j,k=0,\ldots,p$, and
similarly for $\bm{d}^{(2)}_{ijk}$ and $\bm{d}^{(3)}_{ijk}$ along the
other two coordinate directions. These are essentially numerical
approximations to the gradient $\bm{q}=\nabla\bm{u}$, but again we
point out that they might depend implicitly on $\bm{q}$ through the
numerical fluxes (if $C_{22}\ne 0$). Our final semi-discrete
formulation for (\ref{splitdisc1}), (\ref{splitdisc2}) gets the form
\begin{align}
\frac{d\bm{u}_{ijk}}{dt} + \frac{1}{J_{ijk}}\sum_{n=1}^3 \bm{r}_{ijk}^{(n)} &= \bm{S}(\bm{u}_{ijk},\bm{q}_{ijk}) \label{semidisc2-1} \\
\frac{1}{J_{ijk}}\sum_{n=1}^3 \bm{d}_{ijk}^{(n)} &= \bm{q}_{ijk} \label{semidisc2-2}
\end{align}

\section{Temporal discretization and nonlinear solvers}

\subsection{Method of lines and time integration} 

We use various techniques to solve the semi-discrete system of equations
(\ref{semidisc2-1}), (\ref{semidisc2-2}), either by integrating in time or solving
for steady-state solutions. First, we define the vectors $\bm{U},\bm{Q}$ with
all solution components $\bm{u}_{ijk},\bm{q}_{ijk}$, respectively, and write
the system as
\begin{align}
\frac{d\bm{U}}{dt} &= \bm{R}(\bm{U},\bm{Q}), \label{semidiscvector1} \\
\bm{Q} &= \bm{D}(\bm{U},\bm{Q}). \label{semidiscvector2}
\end{align}
This split form can be useful for implicit time-stepping or steady-state
solutions, in particular if the coefficient $C_{22}\ne 0$.  With a standard
Newton's method, this requires the solution of linear systems involving the
matrix
\begin{align}
\bm{K} =
\begin{bmatrix}
\frac{\partial \bm{R}}{\partial \bm{U}} & \frac{\partial \bm{R}}{\partial \bm{Q}} \\
\frac{\partial \bm{D}}{\partial \bm{U}} & \frac{\partial \bm{D}}{\partial \bm{Q}}
\end{bmatrix} \equiv
\begin{bmatrix}
\bm{K}_{11} & \bm{K}_{12} \\
\bm{K}_{21} & \bm{K}_{22}
\end{bmatrix}.
\end{align}
This system solves for both $\bm{U}$ and $\bm{Q}$ but it retains the high
level of sparsity of the method. Also, it allows for non-zero $\bm{K}_{22}$
which can be used to give the scheme several attractive properties
\cite{cockburn08hybrid}. In our examples, we solve these equations using a
standard sparse direct solver \cite{davis06sparse}.

However, in most of our problems we set $C_{22}=0$ to allow for elimination of
the discrete derivatives $\bm{Q}$. Then
$\bm{D}(\bm{U},\bm{Q})=\bm{D}(\bm{U})$, and substituting
(\ref{semidiscvector2}) into (\ref{semidiscvector1}) leads to a reduced
system
\begin{align}
\frac{d\bm{U}}{dt} &= \bm{R}(\bm{U},\bm{D}(\bm{U})) \equiv \bm{F}(\bm{U}). \label{semidiscprimal}
\end{align}
This is clearly the preferred choice for explicit time-stepping, since it is a
regular system of ODEs. In our examples we use a standard fourth-order
explicit Runge-Kutta method. We also use this form for implicit time-stepping
using Diagonally Implicit Runge-Kutta (DIRK) schemes
\cite{alexander77dirk}. In particular, we use the following L-stable,
three-stage, third-order accurate method \cite{alexander77dirk}:
\begin{align}
  \bm{K}_i &= \bm{F}\bigg(\bm{U}_n+\Delta t\sum_{j=1}^s a_{ij}
  \bm{K}_j\bigg),\quad i=1,\ldots,s \label{dirk1} \\ \bm{U}_{n+1} &= \bm{U}_n + \Delta t\sum_{j=1}^s b_j
  \bm{K}_j, \label{dirk2}
\end{align}
with $s=3$ and the coefficients given by the Runge-Kutta tableaux below.
\begin{center}
  \begin{minipage}{.12\textwidth}
    \vspace{11mm}
    \begin{tabular}{c|c}
      $c$ & $A$ \\ \hline
      & $b^T$
    \end{tabular} $=$
  \end{minipage}
  \begin{minipage}{.22\textwidth}
    \vspace{0mm}
    \renewcommand{\arraystretch}{1.1}
    \begin{tabular}{c|ccc}
      $\alpha$ & $\alpha$ & 0 & \ \ 0 \\
      $\tau_2$ & $\tau_2-\alpha$ & $\alpha$ & \ \ 0 \\
      \rule[-5pt]{0pt}{0pt}        $1$ & $b_1$ & $b_2$ & \ \ $\alpha$ \\ \hline
      \rule[12pt]{0pt}{0pt}    & $b_1$ & $b_2$ & $\alpha$
    \end{tabular}
  \end{minipage}
  \hspace{10mm}
  \begin{minipage}{.25\textwidth}
    \vspace{-4mm}
    \begin{align*}
      \alpha &= 0.435866521508459\\
      \tau_2 &= (1+\alpha)/2\\
      b_1 &= -(6\alpha^2-16\alpha+1)/4 \\
      b_2 &= (6\alpha^2-20\alpha+5)/4
    \end{align*}
  \end{minipage}
\end{center}
We also use implicit time-stepping for computing steady-state solutions, by a
sequence of increasing timesteps $\Delta t$ and a final step without the time
derivatives. Since this does not require time-accuracy, we use a standard
backward Euler scheme. We solve the nonlinear systems (\ref{dirk1}) using
Newton's method with preconditioned iterative solvers, as described below.

\subsection{Newton-Krylov solvers}
\label{sec:newtonkrylov}

When Newton's method is applied to the reduced problem (\ref{semidiscprimal}),
it requires the solution of systems of equations of the form
\begin{align}
(\bm{I}-\alpha \Delta t \bm{A}) \Delta \bm{U}^{(i)} = \Delta t \bm{R} (\bm{U}^{(i)},\bm{D}(\bm{U}^{(i)})) \label{dirkeqn}
\end{align}
where $\Delta t$ is the timestep and
\begin{align}
\bm{A} = \frac{d \bm{R}}{d\bm{U}} = \frac{\partial \bm{R}}{\partial \bm{U}} +
\frac{\partial \bm{R}}{\partial \bm{Q}} \frac{\partial \bm{D}}{\partial \bm{U}}
= \bm{K}_{11} + \bm{K}_{12}\bm{K}_{21}. \label{fulldiscprimal}
\end{align}

Forming this matrix $\bm{A}$ has the drawback that for second-order systems,
the product $\bm{K}_{12}\bm{K}_{21}$ is in general much less sparse than the
individual matrices $\bm{K}_{11},\bm{K}_{12},\bm{K}_{21}$. This is expected
due to the repeated differentiation along two different directions, but it
requires special solvers to avoid explicitly forming the denser matrix
$\bm{A}$. This phenomenon is not unique for our method, in fact many other
numerical schemes including finite difference methods and nodal DG methods
suffer from sparsity reduction for second-order systems.

In this work, we use a simple approach to solve the system (\ref{dirkeqn})
without forming the full Jacobian matrix. In a preconditioned Krylov subspace
method, we need to perform two operations: Multiplication of a vector $\bm{p}$
by the matrix $(\bm{I}-\alpha \Delta t \bm{A})$, and approximate solution of
$(\bm{I}-\alpha \Delta t \bm{A})\bm{x}=\bm{b}$ for the preconditioning.  The
matrix-vector product can by computed by keeping the individual matrix in a
separated form and nesting the products:
\begin{align}
(\bm{I}-\alpha \Delta t \bm{A})\bm{p} = \bm{p}-\alpha\Delta t \left(\bm{K}_{11}\bm{p} + \bm{K}_{12} (\bm{K}_{21} \bm{p}) \right). \label{matvecsplit}
\end{align}
This avoids explicitly forming the matrix $\bm{A}$, and the cost per
matrix-vector product is proportional to the number of entries in the matrices
$\bm{K}_{11},\bm{K}_{12},\bm{K}_{21}$.

For preconditioning, we use a sparse block-Jacobi approach which forms an
approximate matrix $\widetilde{\bm{A}}$ that ignores all fill from the product
$\bm{K}_{12}\bm{K}_{21}$ and any inter-element connectivities. In other words,
$\widetilde{\bm{A}}$ is equal to $\bm{A}$ only at the block-diagonal Line-DG
sparsity pattern, and zero everywhere else. This simple preconditioner
requires very little storage (even less than a first-order discretization or
the matrix $\bm{K}_{11}$) and it performs well for time-accurate simulations
with small timesteps. It can certainly be improved upon, for example allowing
higher levels of fill or using block-ILU and $h/p$-multigrid
schemes \cite{persson08newtongmres}.

When solving the linear systems involving the preconditioning matrix
$(\bm{I}-\alpha \Delta t \widetilde{\bm{A}})$, we use a sparse direct
LU-factorization with fill-reducing ordering for each block
\cite{davis06sparse}. This results in some additional fill, but in our 2-D
examples we find that even for $p$ as high as $7$ the number of entries in
$\widetilde{\bm{A}}$ is about the same as in the original sparse matrices
$\bm{K}_{11}$, $\bm{K}_{12}$, and $\bm{K}_{21}$. For 3-D problems, it is more
critical to retain the line-based sparsity in the preconditioner, and a number
of alternatives should be applicable such as low-order approximations
\cite{orszag80spectral}, ADI-iterations \cite{canuto90adi}, and subiterations
\cite{rumsey95subiterations}.

In our implementation, we store all matrices in a general purpose compressed
column storage format \cite{davis06sparse}. We also point out that the matrix
$\bm{K}_{21}$ can be handled very efficiently, since it is a discrete gradient
operator and therefore (a) linear, (b) constant in time, and (c) equal for all
solution components (except possibly at the boundaries for certain boundary
conditions).

\subsection{Re-using Jacobian matrices}
\label{sec:reuse}

In many problems it is more computationally expensive to form the matrix
$\bm{I}-\alpha\Delta t \bm{A}$ than to solve the linear system
(\ref{dirkeqn}). This is especially true for time-accurate integration where
the timesteps $\Delta t$ are relatively small (but still large enough to
motivate the use of implicit solvers). We therefore use a standard technique
for Newton's method that attempts to re-use old Jacobian matrices for many
iterations, until the convergence is too slow as determined by the number of
iterations exceeding a threshold number. In our numerical experiments this is
sufficient to allow for a reuse of the Jacobian for a large number of Newton
steps, DIRK stages, and timesteps at a time.

The linear systems are solved using a preconditioned GMRES method
\cite{saad86gmres}, with the low-cost sparse block-Jacobi preconditioner
$\widetilde{\bm{A}}$ described above.  These equations can be solved with
relatively low accuracy using a small number of GMRES iterations, since we are
using old Jacobian matrices which already limit the potential improvements
from each Newton step. The tolerance in the Newton solver is set well below
the estimated truncation error of the time integrator.

\section{Results}

\subsection{Poisson's equation}

Our first test is Poisson's equation
\begin{align}
-\nabla \cdot (\nabla u) = f(x,y) \label{poisson}
\end{align}
on the unit square domain $\Omega = [0,1]^2$. Dirichlet conditions are imposed
at all the boundaries ($\partial \Omega_D = \partial \Omega$) and we choose
the analytical solution
\begin{align}
  u(x,y) = \mathrm{exp} \left[ \alpha \sin (ax+by) + \beta \cos (cx+dy)
                              \right]    \label{solution}
\end{align}
with numerical parameters $\alpha=0.1, \beta=0.3, a=5.1, b=-6.2, c=4.3,
d=3.4$. We then solve (\ref{poisson}) with Dirichlet boundary conditions
$g_D(x,y)=u(x,y)|_{\partial \Omega_D}$. The source term, $f(x,y)$, is obtained
by analytical differentiation of (\ref{solution}).

\begin{figure}[t]
\begin{center}
\begin{minipage}{.3\textwidth}
  \begin{center}
    \includegraphics[width=.8\textwidth]{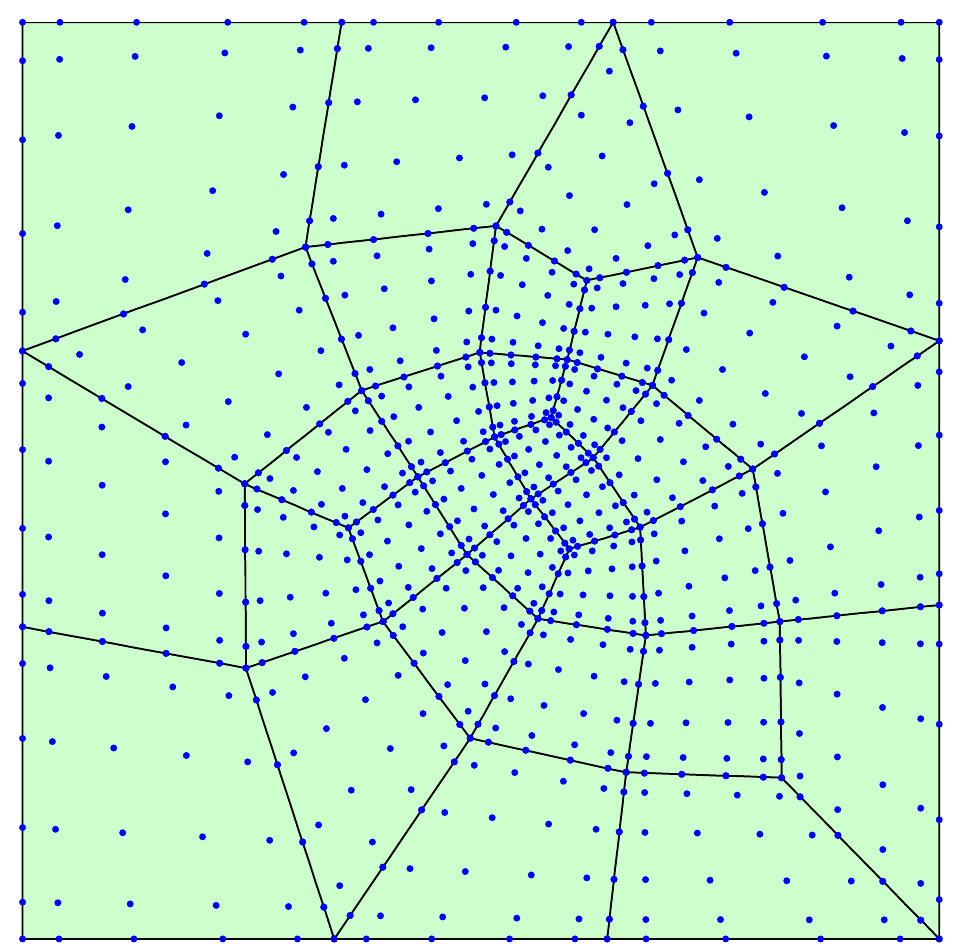} \\
    Coarsest mesh, $p=5$
  \end{center}
\end{minipage}
\begin{minipage}{.3\textwidth}
  \begin{center}
    \includegraphics[width=.8\textwidth]{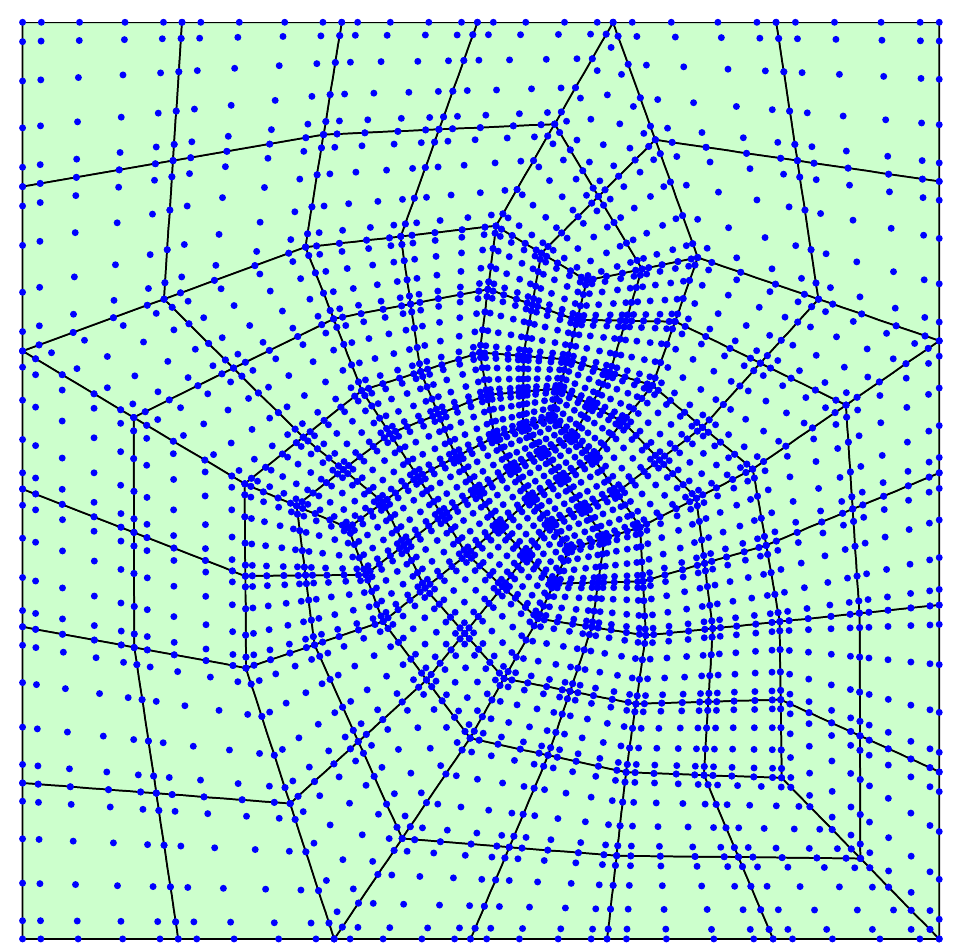} \\
    One refinement, $p=5$
  \end{center}
\end{minipage}
\begin{minipage}{.3\textwidth}
  \begin{center}
    \includegraphics[width=.8\textwidth]{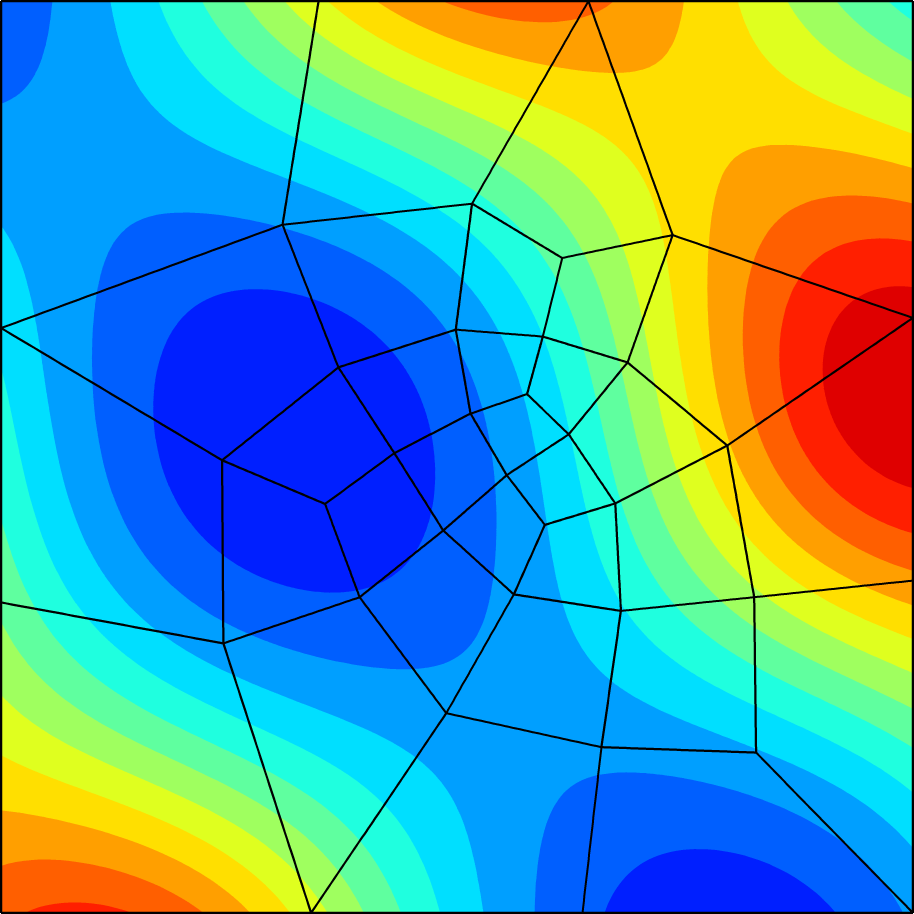} \\
    Solution $u(x,y)$
  \end{center}
\end{minipage}
\end{center}
\caption{The Poisson test problem (\ref{poisson}). The left figure shows
the coarse unstructured quadrilateral mesh, which is uniformly refined
repeatedly, and the solution nodes for $p=5$. The right figure shows the
solution as colored contours.}
\label{poitest}
\end{figure}

We discretize $\Omega$ using an unstructured mesh of quadrilateral elements,
see Figure~\ref{poitest} (left). We solve the split system
(\ref{semidiscvector1}), (\ref{semidiscvector2}) using a direct sparse solver,
for polynomial degrees $p=1,\ldots,7$. We consider two sets of parameters: the
minimal dissipation scheme with $C_{11}=C_{22}=0$, which has the benefit that
it allows for elimination of the gradients $\bm{q}$, and the slightly
over-stabilized scheme $C_{11}=C_{22}=1/20$, which may provide a higher-order
of convergence for $\bm{q}$ \cite{cockburn08hybrid}.

\begin{table}[h]
    \begin{footnotesize}
  \begin{center}
    \textbf{Error in solution $|u_h-u|_\infty$, $C_{11}=C_{22}=0$, $\bm{C}_{12}=\bm{n}S_i^\pm/2$} \\
    \begin{tabular}{r|l@{\ }r@{\ \ }|l@{\ }r@{\ \ }|l@{\ }r@{\ \ }|l@{\ }r@{\ \ }|l@{\ }r@{\ \ }|l@{\ }r@{\ \ }|l@{\ }r}
      &
      \multicolumn{2}{c}{$p=1$} &  \multicolumn{2}{c}{$p=2$} & \multicolumn{2}{c}{$p=3$} & \multicolumn{2}{c}{$p=4$} & \multicolumn{2}{c}{$p=5$} & \multicolumn{2}{c}{$p=6$} & \multicolumn{2}{c}{$p=7$} \\
$n$
& Error & Rate & Error & Rate & Error & Rate & Error & Rate & Error & Rate & Error & Rate & Error & Rate \\ \hline
1 & $4.6 \cdot 10^{-2}$ &     & $3.7 \cdot 10^{-3}$ &     & $4.1 \cdot 10^{-4}$ &     & $1.0 \cdot 10^{-4}$ &     & $8.9 \cdot 10^{-6}$ &     & $2.1 \cdot 10^{-6}$ &     & $3.7 \cdot 10^{-7}$ &     \\
2 & $1.1 \cdot 10^{-2}$ & 2.0 & $4.3 \cdot 10^{-4}$ & 3.1 & $2.8 \cdot 10^{-5}$ & 3.9 & $2.5 \cdot 10^{-6}$ & 5.4 & $1.3 \cdot 10^{-7}$ & 6.0 & $1.5 \cdot 10^{-8}$ & 7.1 & $9.4 \cdot 10^{-10}$ & 8.6 \\
4 & $2.7 \cdot 10^{-3}$ & 2.1 & $5.1 \cdot 10^{-5}$ & 3.1 & $1.7 \cdot 10^{-6}$ & 4.0 & $6.0 \cdot 10^{-8}$ & 5.4 & $2.1 \cdot 10^{-9}$ & 6.0 & $8.1 \cdot 10^{-11}$ & 7.5 & $2.7 \cdot 10^{-12}$ & 8.4 \\
8 & $6.6 \cdot 10^{-4}$ & 2.0 & $6.2 \cdot 10^{-6}$ & 3.0 & $1.0 \cdot 10^{-7}$ & 4.0 & $1.8 \cdot 10^{-9}$ & 5.1 & $3.0 \cdot 10^{-11}$ & 6.1 & $5.4 \cdot 10^{-13}$ & 7.2 & $9.1 \cdot 10^{-14}$ & * \\
    \end{tabular}
\\ \ \\ \ \\
    \textbf{Error in gradient $|\bm{q}_h-\nabla u|_\infty$, $C_{11}=C_{22}=0$, $\bm{C}_{12}=\bm{n}S_i^\pm/2$} \\
    \begin{tabular}{r|l@{\ }r@{\ \ }|l@{\ }r@{\ \ }|l@{\ }r@{\ \ }|l@{\ }r@{\ \ }|l@{\ }r@{\ \ }|l@{\ }r@{\ \ }|l@{\ }r}
      &
      \multicolumn{2}{c}{$p=1$} &  \multicolumn{2}{c}{$p=2$} & \multicolumn{2}{c}{$p=3$} & \multicolumn{2}{c}{$p=4$} & \multicolumn{2}{c}{$p=5$} & \multicolumn{2}{c}{$p=6$} & \multicolumn{2}{c}{$p=7$} \\
$n$
& Error & Rate & Error & Rate & Error & Rate & Error & Rate & Error & Rate & Error & Rate & Error & Rate \\ \hline
1 & $4.6 \cdot 10^{-2}$ &     & $3.7 \cdot 10^{-3}$ &     & $4.1 \cdot 10^{-4}$ &     & $1.0 \cdot 10^{-4}$ &     & $8.9 \cdot 10^{-6}$ &     & $2.1 \cdot 10^{-6}$ &     & $3.7 \cdot 10^{-7}$ &     \\
2 & $1.1 \cdot 10^{-2}$ & 1.1 & $4.3 \cdot 10^{-4}$ & 2.2 & $2.8 \cdot 10^{-5}$ & 3.2 & $2.5 \cdot 10^{-6}$ & 4.4 & $1.3 \cdot 10^{-7}$ & 5.0 & $1.5 \cdot 10^{-8}$ & 6.4 & $9.4 \cdot 10^{-10}$ & 7.3 \\
4 & $2.7 \cdot 10^{-3}$ & 1.0 & $5.1 \cdot 10^{-5}$ & 2.1 & $1.7 \cdot 10^{-6}$ & 3.0 & $6.0 \cdot 10^{-8}$ & 4.2 & $2.1 \cdot 10^{-9}$ & 5.0 & $8.1 \cdot 10^{-11}$ & 6.4 & $2.7 \cdot 10^{-12}$ & 6.9 \\
8 & $6.6 \cdot 10^{-4}$ & 1.0 & $6.2 \cdot 10^{-6}$ & 2.0 & $1.0 \cdot 10^{-7}$ & 2.9 & $1.8 \cdot 10^{-9}$ & 4.0 & $3.0 \cdot 10^{-11}$ & 5.0 & $5.4 \cdot 10^{-13}$ & * & $9.1 \cdot 10^{-14}$ & * \\
    \end{tabular}
\\ \ \\ \ \\
    \textbf{Error in solution $|u_h-u|_\infty$, $C_{11}=C_{22}=1/20$, $\bm{C}_{12}=\bm{n}S_i^\pm/2$} \\
    \begin{tabular}{r|l@{\ }r@{\ \ }|l@{\ }r@{\ \ }|l@{\ }r@{\ \ }|l@{\ }r@{\ \ }|l@{\ }r@{\ \ }|l@{\ }r@{\ \ }|l@{\ }r}
      &
      \multicolumn{2}{c}{$p=1$} &  \multicolumn{2}{c}{$p=2$} & \multicolumn{2}{c}{$p=3$} & \multicolumn{2}{c}{$p=4$} & \multicolumn{2}{c}{$p=5$} & \multicolumn{2}{c}{$p=6$} & \multicolumn{2}{c}{$p=7$} \\
$n$
& Error & Rate & Error & Rate & Error & Rate & Error & Rate & Error & Rate & Error & Rate & Error & Rate \\ \hline
1 & $4.5 \cdot 10^{-2}$ &     & $3.9 \cdot 10^{-3}$ &     & $4.9 \cdot 10^{-4}$ &     & $1.2 \cdot 10^{-4}$ &     & $9.3 \cdot 10^{-6}$ &     & $2.5 \cdot 10^{-6}$ &     & $3.8 \cdot 10^{-7}$ &     \\
2 & $1.0 \cdot 10^{-2}$ & 2.1 & $4.9 \cdot 10^{-4}$ & 3.0 & $3.7 \cdot 10^{-5}$ & 3.7 & $3.2 \cdot 10^{-6}$ & 5.2 & $1.8 \cdot 10^{-7}$ & 5.7 & $2.2 \cdot 10^{-8}$ & 6.9 & $1.1 \cdot 10^{-9}$ & 8.4 \\
4 & $2.4 \cdot 10^{-3}$ & 2.1 & $5.8 \cdot 10^{-5}$ & 3.1 & $2.3 \cdot 10^{-6}$ & 4.0 & $7.7 \cdot 10^{-8}$ & 5.4 & $3.4 \cdot 10^{-9}$ & 5.7 & $1.2 \cdot 10^{-10}$ & 7.4 & $4.5 \cdot 10^{-12}$ & 8.0 \\
8 & $5.9 \cdot 10^{-4}$ & 2.0 & $7.2 \cdot 10^{-6}$ & 3.0 & $1.3 \cdot 10^{-7}$ & 4.1 & $2.2 \cdot 10^{-9}$ & 5.1 & $4.5 \cdot 10^{-11}$ & 6.2 & $8.2 \cdot 10^{-13}$ & 7.2 & $1.4 \cdot 10^{-13}$ & * \\
    \end{tabular}
\\ \ \\ \ \\
    \textbf{Error in gradient $|\bm{q}_h-\nabla u|_\infty$, $C_{11}=C_{22}=1/20$, $\bm{C}_{12}=\bm{n}S_i^\pm/2$} \\
    \begin{tabular}{r|l@{\ }r@{\ \ }|l@{\ }r@{\ \ }|l@{\ }r@{\ \ }|l@{\ }r@{\ \ }|l@{\ }r@{\ \ }|l@{\ }r@{\ \ }|l@{\ }r}
      &
      \multicolumn{2}{c}{$p=1$} &  \multicolumn{2}{c}{$p=2$} & \multicolumn{2}{c}{$p=3$} & \multicolumn{2}{c}{$p=4$} & \multicolumn{2}{c}{$p=5$} & \multicolumn{2}{c}{$p=6$} & \multicolumn{2}{c}{$p=7$} \\
$n$
& Error & Rate & Error & Rate & Error & Rate & Error & Rate & Error & Rate & Error & Rate & Error & Rate \\ \hline
1 & $4.5 \cdot 10^{-2}$ &     & $3.9 \cdot 10^{-3}$ &     & $4.9 \cdot 10^{-4}$ &     & $1.2 \cdot 10^{-4}$ &     & $9.3 \cdot 10^{-6}$ &     & $2.5 \cdot 10^{-6}$ &     & $3.8 \cdot 10^{-7}$ &     \\
2 & $1.0 \cdot 10^{-2}$ & 1.8 & $4.9 \cdot 10^{-4}$ & 2.5 & $3.7 \cdot 10^{-5}$ & 3.7 & $3.2 \cdot 10^{-6}$ & 4.6 & $1.8 \cdot 10^{-7}$ & 5.6 & $2.2 \cdot 10^{-8}$ & 6.6 & $1.1 \cdot 10^{-9}$ & 7.7 \\
4 & $2.4 \cdot 10^{-3}$ & 1.7 & $5.8 \cdot 10^{-5}$ & 2.5 & $2.3 \cdot 10^{-6}$ & 3.7 & $7.7 \cdot 10^{-8}$ & 4.6 & $3.4 \cdot 10^{-9}$ & 5.6 & $1.2 \cdot 10^{-10}$ & 6.6 & $4.5 \cdot 10^{-12}$ & 7.6 \\
8 & $5.9 \cdot 10^{-4}$ & 1.8 & $7.2 \cdot 10^{-6}$ & 2.6 & $1.3 \cdot 10^{-7}$ & 3.7 & $2.2 \cdot 10^{-9}$ & 4.5 & $4.5 \cdot 10^{-11}$ & 5.7 & $8.2 \cdot 10^{-13}$ & * & $1.4 \cdot 10^{-13}$ & * \\
    \end{tabular}
  \end{center}
    \end{footnotesize}
    \caption{Convergence of $u$ and $\bm{q}$ for the Poisson problem, with $C_{12} = \bm{n}S^\pm_i /2$
    and $C_{11}=C_{22}=0$ (top two tables), $C_{11}=C_{22}=1/20$ (bottom two tables). For the first case we observe approximate rates of $p+1$ for $u$ and $p$ for $\bm{q}$, while the nonzero $C_{11},C_{22}$ case appears to give a significantly higher rate for $\bm{q}$. A star symbol (*) indicates that the error is dominated by floating point rounding errors rather than the truncation error of the scheme.}
  \label{poitab1}
\end{table}

The resulting infinity norm errors and rates of convergence are shown in
table~\ref{poitab1}, for both the solution $u$ and the gradients $\bm{q}$ and
the two parameter cases. For the minimal dissipation scheme (top two tables),
we observe the expected orders of convergence $p+1$ and $p$ for $u$ and
$\bm{q}$, respectively.  For the stabilized scheme (bottom two tables), we
obtain a somewhat higher order for the $\bm{q}$ variables, which could be used
as part of a postprocessing step to further increase the order of convergence
for the solution $u$ \cite{peraire09hybrid}.

\subsection{Euler vortex}

Next we consider the compressible Euler and Navier-Stokes equations, which we
write in the form:
\begin{align}
\frac{\partial \rho}{\partial t}  + \frac{\partial}{\partial x_i}
(\rho u_i) &= 0, \label{ns1} \\
\frac{\partial}{\partial t} (\rho u_i) +
\frac{\partial}{\partial x_i} (\rho u_i u_j+ p)  &=
+\frac{\partial \tau_{ij}}{\partial x_j}
\quad\text{for }i=1,2,3, \label{ns2} \\
\frac{\partial}{\partial t} (\rho E) +
\frac{\partial}{\partial x_i} \left(u_j(\rho E+p)\right) &=
-\frac{\partial q_j}{\partial x_j}
+\frac{\partial}{\partial x_j}(u_j\tau_{ij}), \label{ns3}
\end{align}
where $\rho$ is the fluid density, $u_1,u_2,u_3$ are the velocity
components, and $E$ is the total energy. The viscous stress tensor and
heat flux are given by
\begin{align}
\tau_{ij} = \mu
\left( \frac{\partial u_i}{\partial x_j} +
\frac{\partial u_j}{\partial x_i} -\frac23
\frac{\partial u_k}{\partial x_j} \delta_{ij} \right)
\qquad \text{ and } \qquad
q_j = -\frac{\mu}{\mathrm{Pr}} \frac{\partial}{\partial x_j}
\left( E+\frac{p}{\rho} -\frac12 u_k u_k \right).
\end{align}
Here, $\mu$ is the viscosity coefficient and $\mathrm{Pr = 0.72}$ is
the Prandtl number which we assume to be constant. For an ideal gas,
the pressure $p$ has the form
\begin{align}
p=(\gamma-1)\rho \left( E - \frac12 u_k u_k\right),
\end{align}
where $\gamma$ is the adiabatic gas constant.

Our first model problem is the inviscid flow of a compressible vortex in a
rectangular domain \cite{erlebacher97vortex}. The vortex is initially centered
at $(x_0,y_0)$ and is moving with the free-stream at an angle $\theta$ with
respect to the $x$-axis. The analytic solution at $(x,y,t)$ is given by
\begin{align}
u &= u_\infty \left(\cos\theta - \frac{\epsilon ((y-y_0)-\bar{v}t)}{2\pi r_c}
                    \exp(f/2) \right), &
\rho &= \rho_\infty \left(1 - \frac{\epsilon^2(\gamma-1)M_\infty^2}{8\pi^2}
                    \exp (f) \right)^\frac{1}{\gamma-1}, \\
v &= u_\infty \left(\sin\theta + \frac{\epsilon ((x-x_0)-\bar{u}t)}{2\pi r_c}
                    \exp (f/2) \right), & 
p    &= p_\infty \left(1 - \frac{\epsilon^2(\gamma-1)M_\infty^2}{8\pi^2}
                    \exp (f) \right)^\frac{\gamma}{\gamma-1},
\end{align}
where $f(x,y,t) =
(1-((x-x_0)-\bar{u}t)^2-((y-y_0)-\bar{v}t)^2)/r_c^2$, $M_\infty$ is
the Mach number, $\gamma=c_p/c_v=1.4$, and $u_\infty$, $p_\infty$,
$\rho_\infty$ are free-stream velocity, pressure, and density. The
Cartesian components of the free-stream velocity are
$\bar{u}=u_\infty\cos\theta$ and $\bar{v}=u_\infty\sin\theta$. The
parameter $\epsilon$ measures the strength of the vortex and $r_c$ is
its size.

We use a domain of size 20-by-15, with the vortex initially centered at
$(x_0,y_0)=(5,5)$ with respect to the lower-left corner. The Mach number is
$M_\infty=0.5$, the angle $\theta=\arctan 1/2$, and the vortex has the
parameters $\epsilon=0.3$ and $r_c=1.5$. We use characteristic boundary
conditions and integrate until time $t_0=\sqrt{10^2+5^2}/10$, when the vortex
has moved a relative distance of $(1,1/2)$.

\begin{figure}[t]
\begin{minipage}{.49\textwidth}
  \begin{center}
    \includegraphics[width=.8\textwidth]{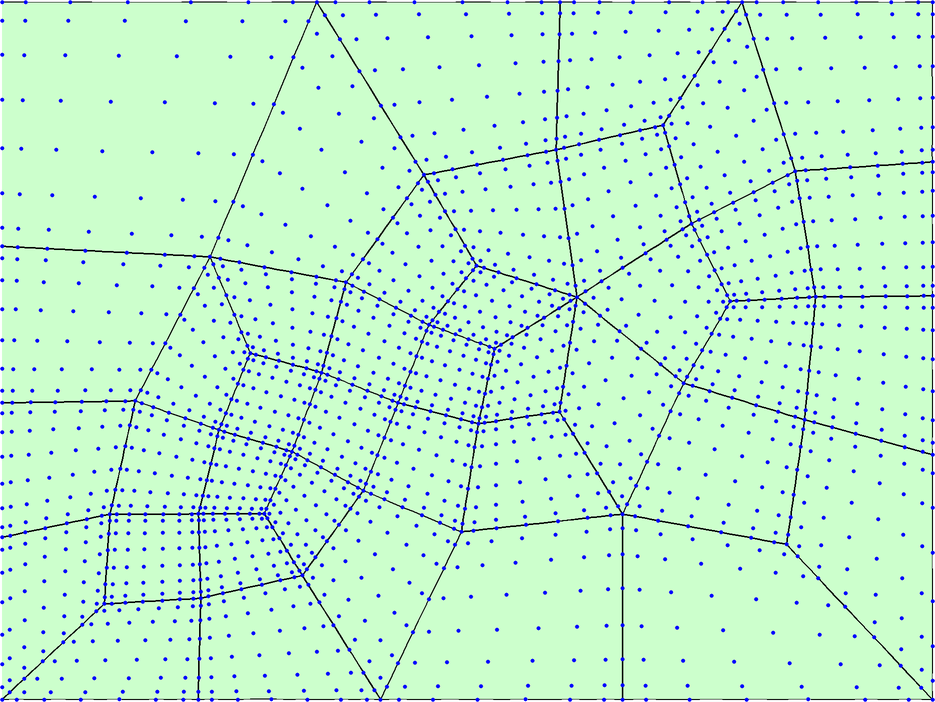} \\
    Coarsest mesh, with degree $p=7$
  \end{center}
\end{minipage}
\begin{minipage}{.49\textwidth}
  \begin{center}
    \includegraphics[width=.8\textwidth]{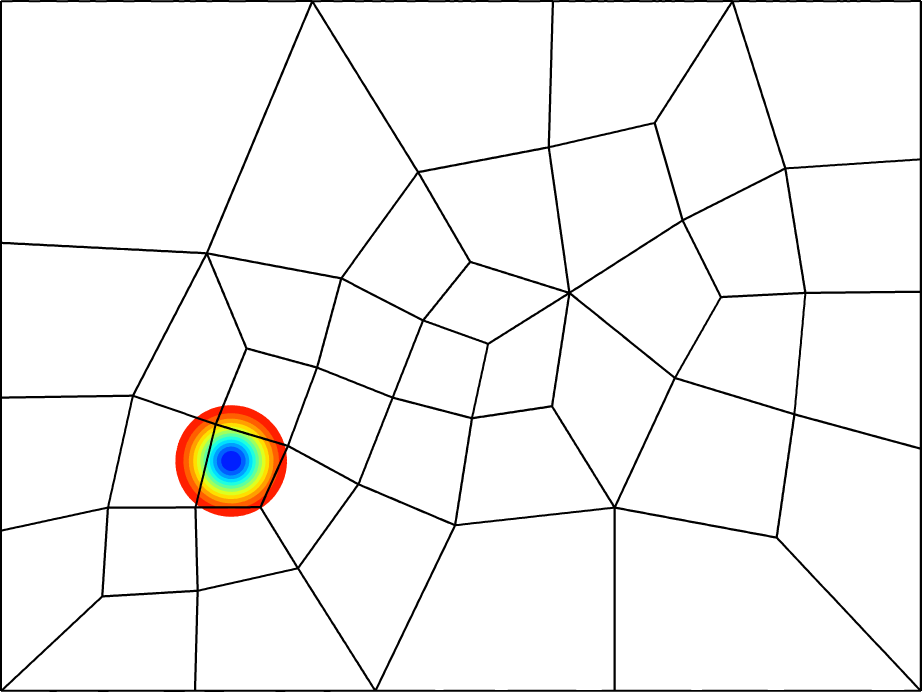} \\
    Solution (density)
  \end{center}
\end{minipage} \\
\begin{center}
  \includegraphics[width=.7\textwidth]{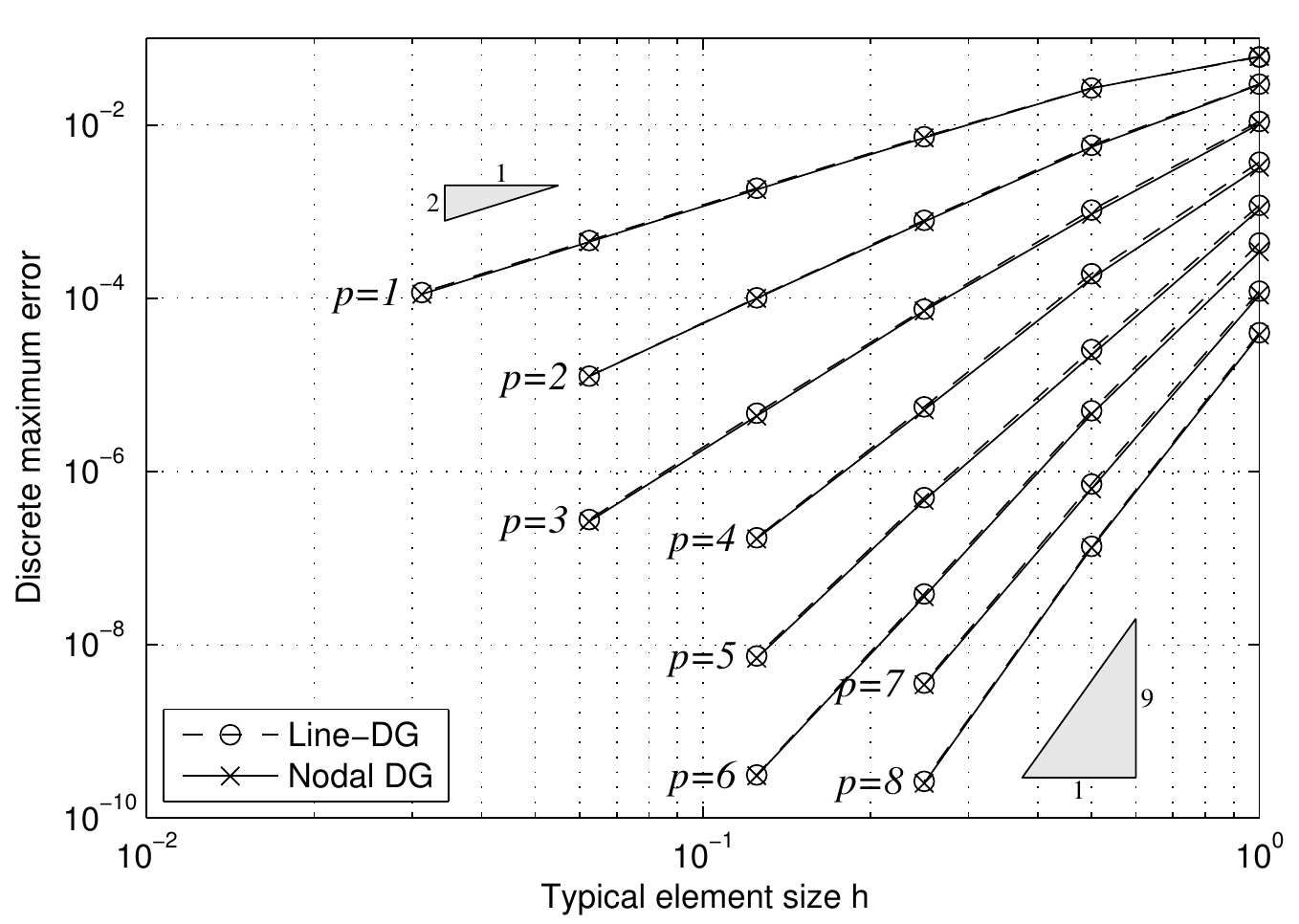}
\end{center}
\caption{Convergence test for an Euler vortex test problem using the Line-DG
  method and the nodal DG method. The results show optimal order of
  convergence $\mathcal{O}(h^{p+1})$ with very small differences between the
  two methods.}
\label{conv}
\end{figure}

We write the Euler equations as a first-order system of conservation laws
(\ref{conslaw}), in the conserved variables $(\rho,\rho u,\rho v,\rho E)$.
The scheme (\ref{semidisc}) is implemented in a straight-forward way, and we
use Roe's method for the numerical fluxes (\ref{numflux1}) \cite{roe}. The
time-integration is done explicitly with the form (\ref{semidiscprimal}) using
the RK4 solver and a timestep $\Delta t$ small enough so that all truncation
errors are dominated by the spatial discretization. We start from a coarse
unstructured quadrilateral mesh (figure~\ref{conv}, top left), which we refine
uniformly a number of times, and we use polynomial degrees $p$ ranging between
1 and 8. The top right plot also shows the density field for a sample
solution.

In the bottom plot of figure~\ref{conv}, we graph the maximum errors
(discretely at the solution nodes) for all simulation cases, both for the
Line-DG method and the standard nodal DG method. The results clearly show the
optimal order of convergence $\mathcal{O}(h^{p+1})$ for element size $h$ for
both methods, and that the Line-DG errors are in all cases very close to those
of the nodal DG method.

\subsection{Inviscid flow over a cylinder}

\begin{figure}[t]
\begin{center}
\begin{minipage}{.34\textwidth}
  \begin{center}
    \includegraphics[width=.9\textwidth]{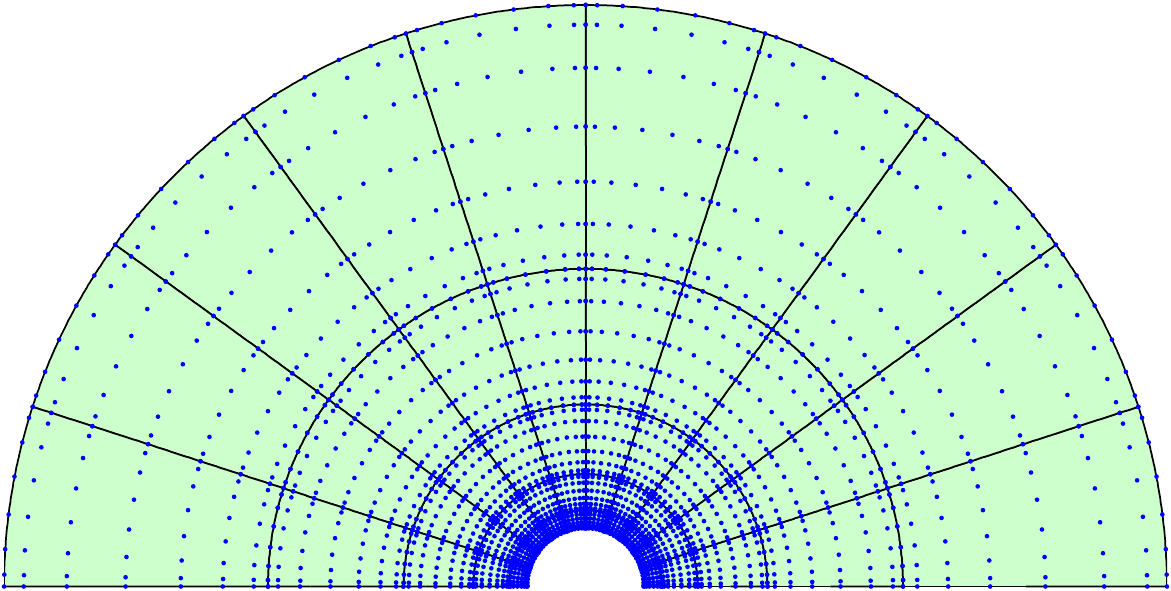} \\
    Coarsest mesh, $p=7$ \\ \ \\
    \includegraphics[width=.9\textwidth]{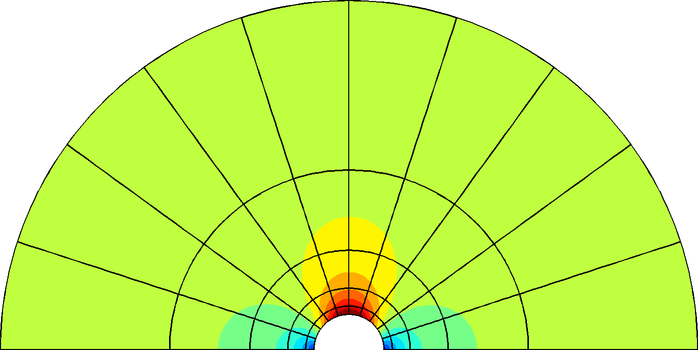} \\ \ \\
    \includegraphics[width=.8\textwidth]{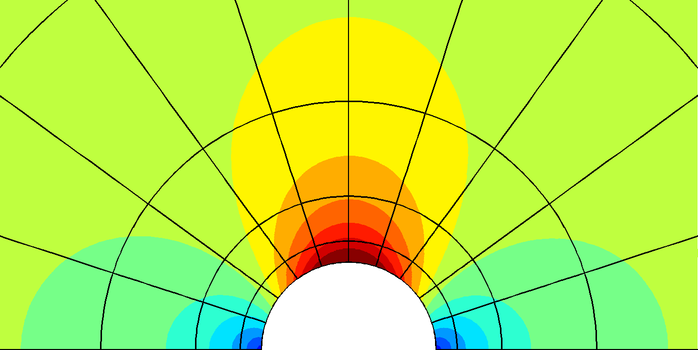} \\ \ \\
    \includegraphics[width=.9\textwidth]{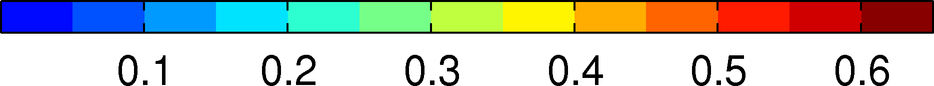} \\
    Solution, Mach number
  \end{center}
\end{minipage} \ \ %
\begin{minipage}{.6\textwidth}
  \begin{center}
    \includegraphics[width=\textwidth]{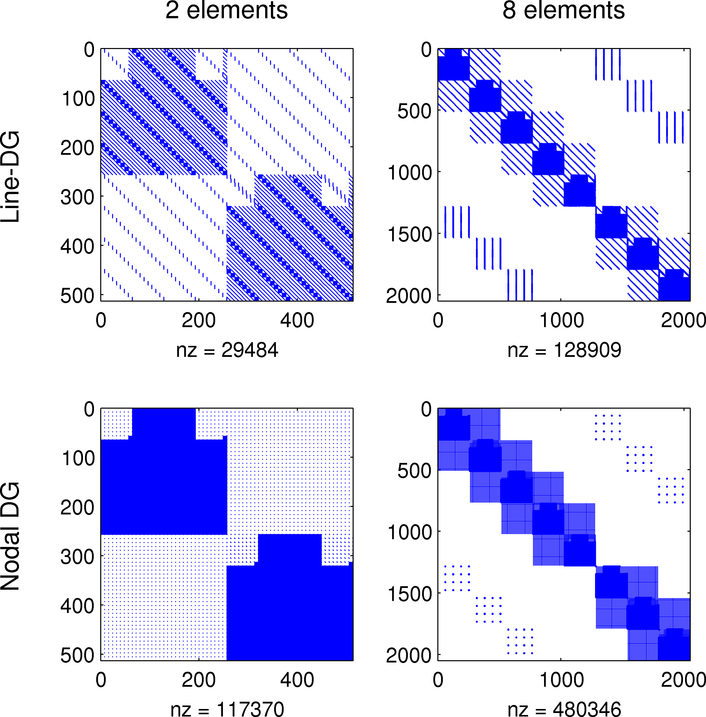} \\
  \end{center}
\end{minipage}
\end{center}
\caption{Inviscid flow over a cylinder. The plots show the coarsest grid used in the convergence study and the nodes for polynomial degree $p=7$ (top left), the corresponding solution as Mach number color plot (bottom left, with zoom-in), and a sparsity plot of the Jacobian matrices for both Line-DG and nodal DG (right).}
\label{pltcylspy}
\end{figure}

Next we study a problem with a steady-state solution and curved boundaries,
and solve the Euler equations for the inviscid flow over a half-cylinder with
radius 1 at a Mach number of 0.3. Structured quadrilateral meshes are used,
with strong element size grading to better resolve the region close to the
cylinder (see figure~\ref{pltcylspy}, top left). The outer domain boundary is
a half-cylinder with radius 10, where characteristic boundary conditions are
imposed. Standard slip wall/symmetry conditions are used at the cylinder and
at the symmetry plane.

The steady-state solutions are found using a fully consistent Newton method,
applied directly to the equations (\ref{semidisc}), with the linear systems
solved using a direct sparse solver \cite{davis06sparse}. Starting the
iterations from an approximate analytical solution, derived from a potential
flow approximation, the solver converges to machine precision in 4 to 6
iterations. The solution is shown in the bottom left of
figure~\ref{pltcylspy}, and the figures to the right show portions of the
Jacobian matrices for both the Line-DG and the nodal DG method. This
illustrates again the reduced sparsity of the Line-DG scheme, with about a
factor of 4 fewer entries than nodal DG already in two space dimensions.

\begin{figure}[h]
\begin{minipage}{.49\textwidth}
  \begin{center}
    \includegraphics[width=.99\textwidth]{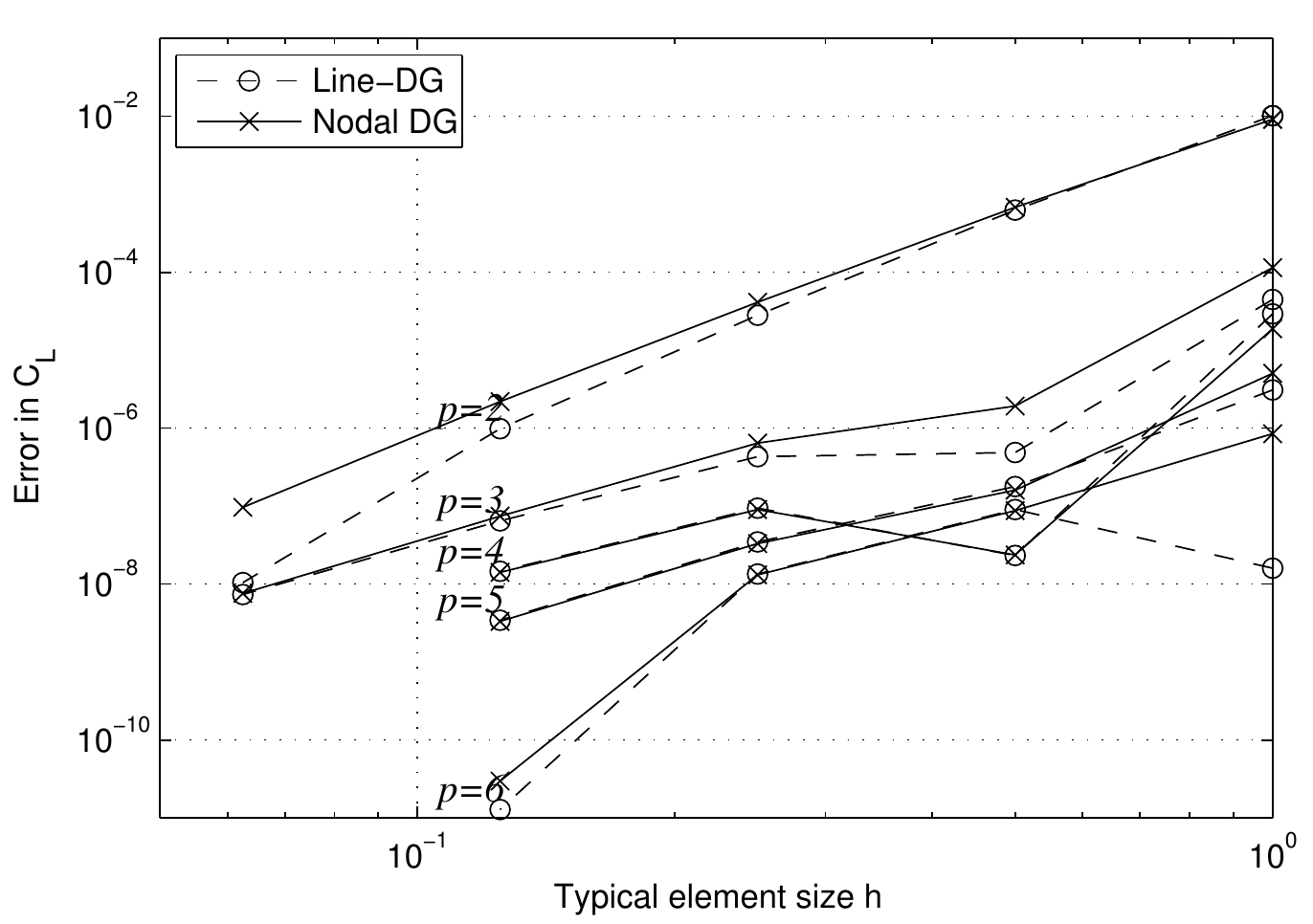} \\
  \end{center}
\end{minipage} \hfill
\begin{minipage}{.49\textwidth}
  \begin{center}
    \includegraphics[width=.99\textwidth]{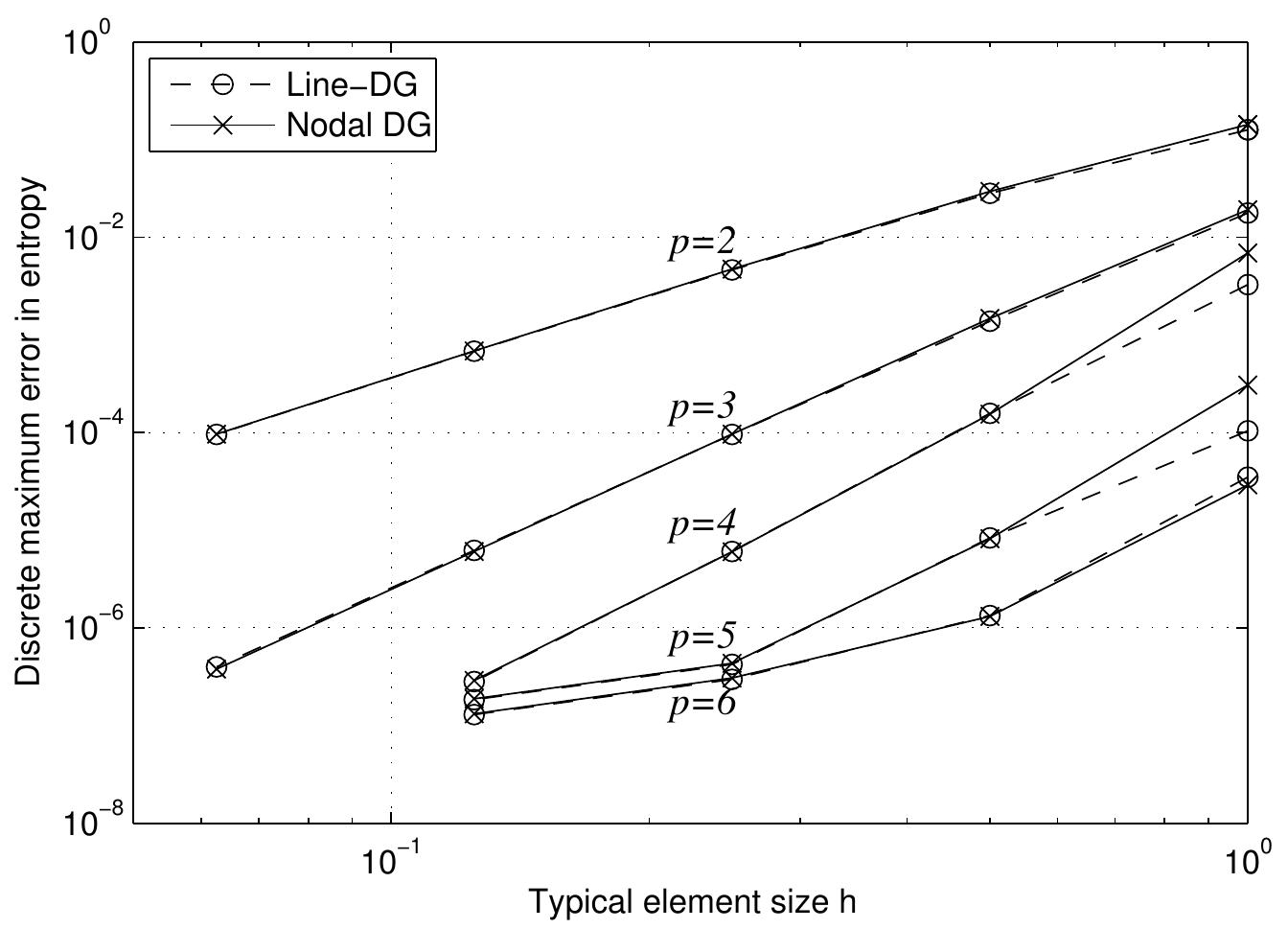} \\
  \end{center}
\end{minipage} \\
\caption{The convergence of the lift coefficient $C_L$ (left) and the
entropy difference (right) for the inviscid flow over cylinder problem. The
plots show a series of results for varying polynomial degrees and number of
refinements, for the two methods Line-DG and nodal DG.}
\label{eulercylconv}
\end{figure}

To evaluate the accuracy and convergence of the scheme, in
figure~\ref{eulercylconv} we plot the errors in the lift coefficient $C_L$
(left) and the maximum errors in the entropy (right). These plots again
confirm the convergence of the schemes as well as the minor differences in
error between the Line-DG and the nodal DG schemes.

\subsection{Laminar flow around airfoil}

\begin{figure}[t!]
  \begin{center}
  \begin{minipage}{.47\textwidth}
    \includegraphics[width=\textwidth]{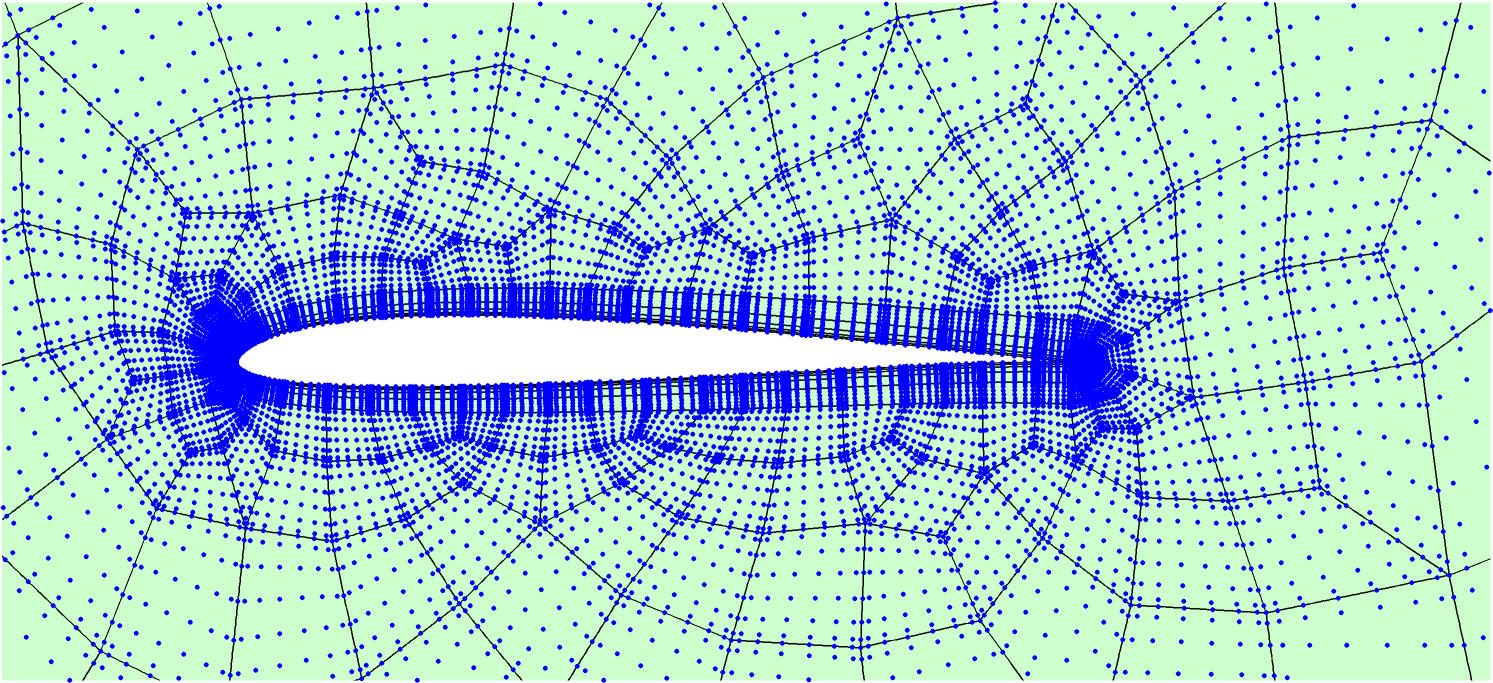}
  \end{minipage} \ \ %
  \begin{minipage}{.47\textwidth}
    \includegraphics[width=\textwidth]{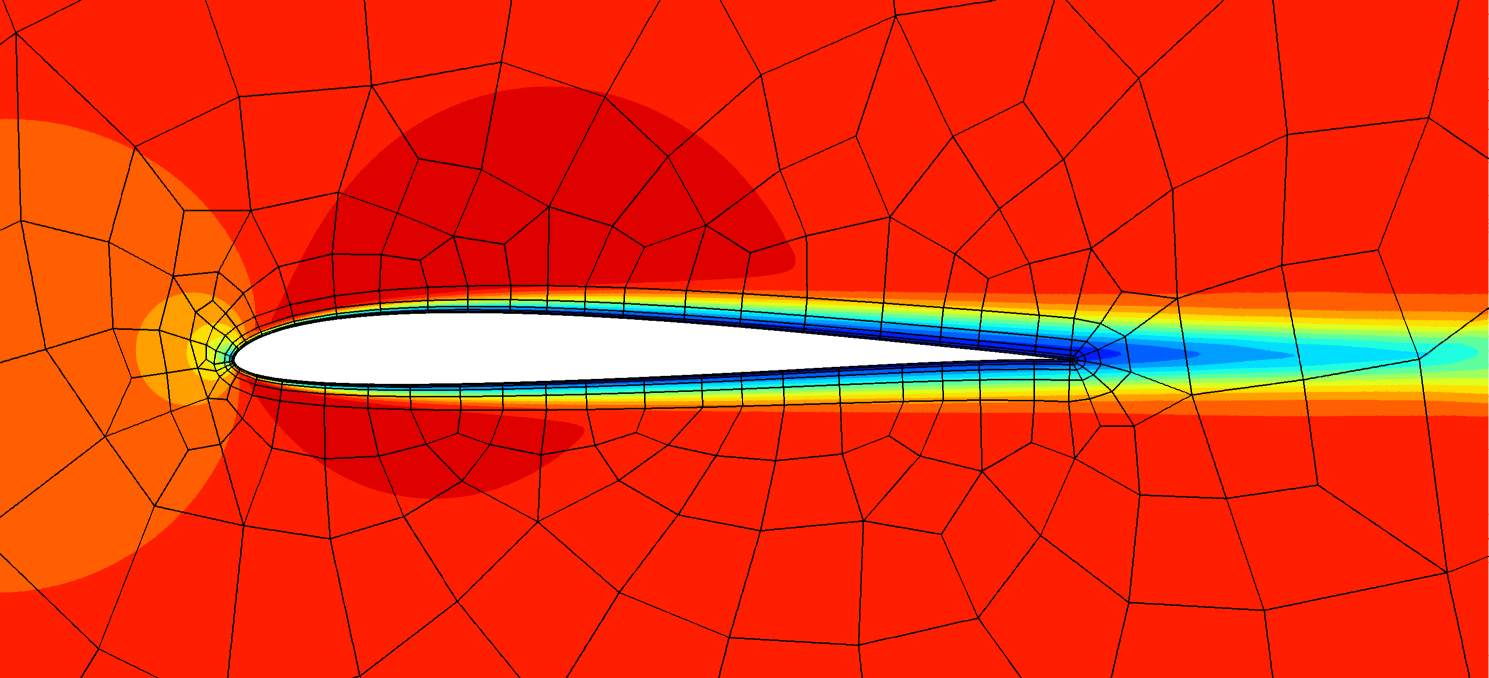}
  \end{minipage} \\ \ \\
  \begin{minipage}{.42\textwidth}
    \includegraphics[width=\textwidth]{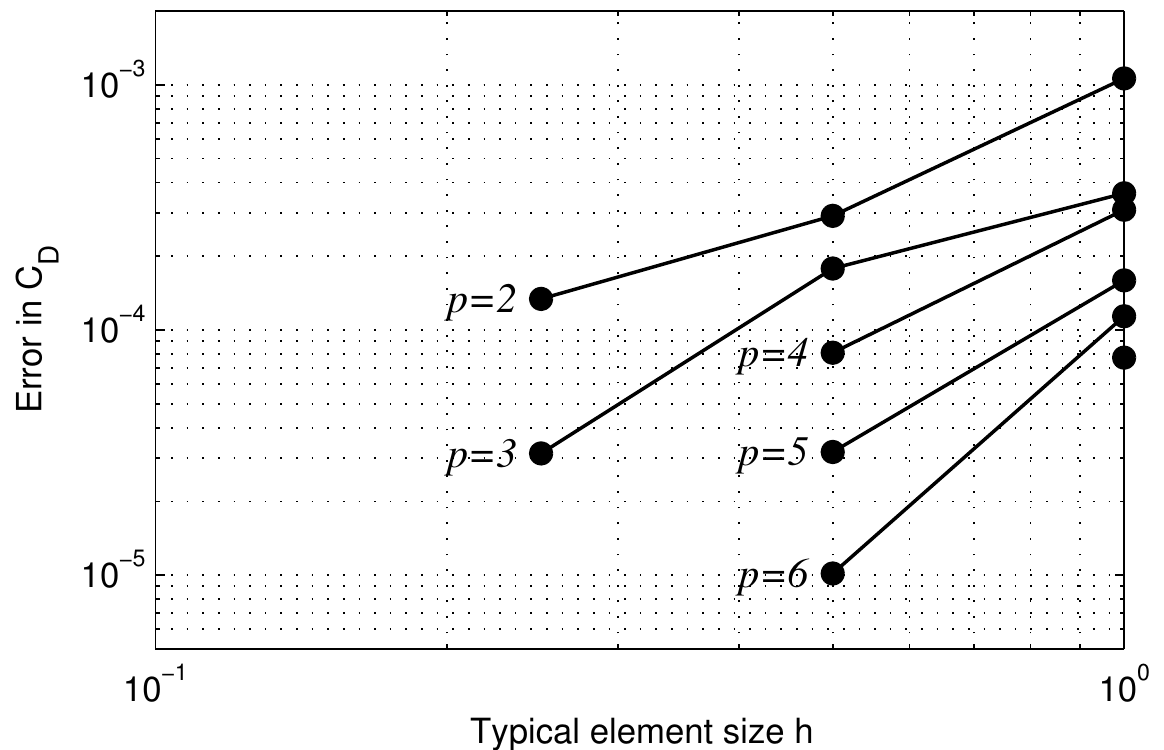}
  \end{minipage} \ \ \ \ \ \ \ %
  \begin{minipage}{.42\textwidth}
    \includegraphics[width=\textwidth]{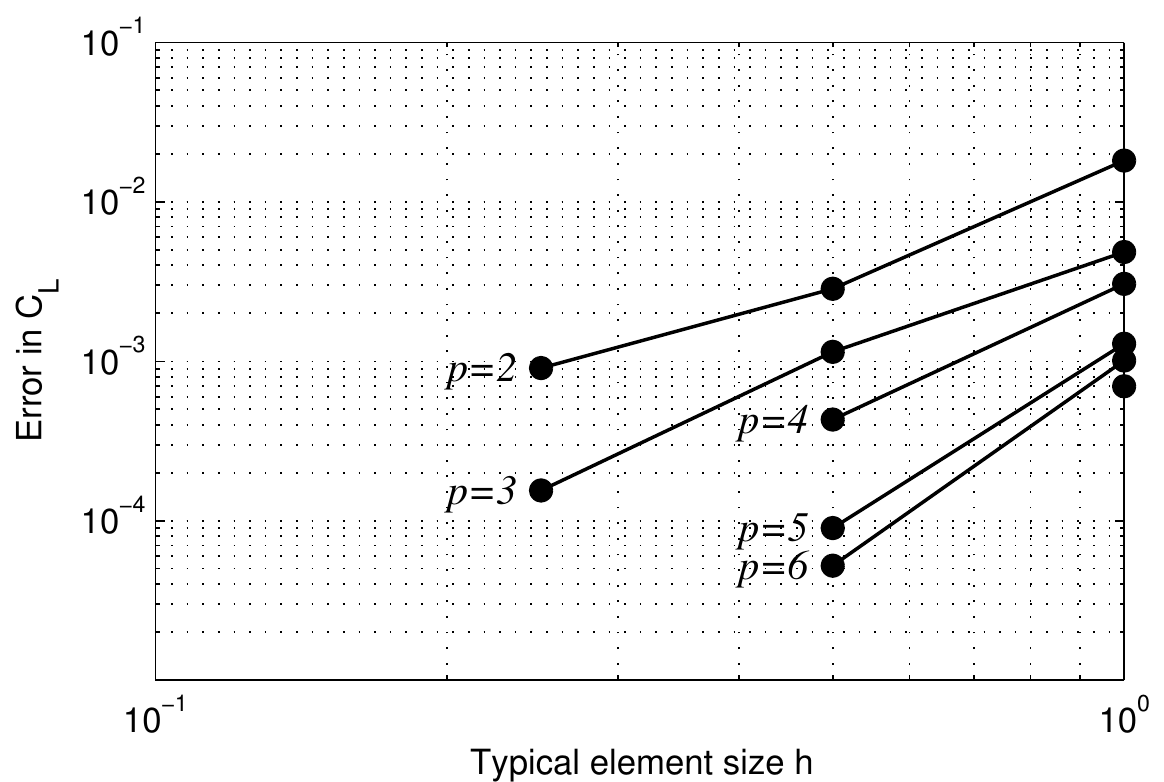}
  \end{minipage}
  \end{center}
  \caption{Stationary flow around an SD7003 airfoil (top left: mesh, top
    right: Mach number), computed with the Line-DG method with $p=7$, at
    free-stream Mach 0.2, zero angle of attack, and Reynolds number 5,000. The
    bottom plots show the convergence of $C_D$ and $C_L$, for a range of
    polynomial degrees and with 0, 1, or 2 uniform mesh refinements.}
  \label{sdfoil}
\end{figure}

An example of a steady-state viscous computation is shown in
figure~\ref{sdfoil}. The compressible Navier-Stokes equations are solved at
Mach 0.2 and Reynolds number 5000, for a flow around an SD7003 airfoil. The
quadrilateral mesh is fully unstructured except for a structured graded
boundary layer region, with a total of 461 elements for the coarse mesh, and
1844 and 7376 elements for the once and twice refined meshes, respectively.
With approximating polynomials of degree $p=6$, this gives a total number of
high-order nodes of 22,589 for the coarse mesh and 90,356 for the first
refinement.

We find the steady-state solution with 10 digits of accuracy in the residual
using a consistent Newton's method, with pseudo-timestepping for
regularization. A solution is shown in figure~\ref{sdfoil} (top right), for
the coarse mesh with $p=7$. In the bottom plots, we show the convergence of
the drag and the lift coefficients, for a range of polynomial degrees
$p$. While it is hard to asses the exact order of convergence from these
numbers, it is clear that our scheme provides a high order of convergence even
for these difficult derivative-based quantities.

\subsection{Transient flow around airfoil at Re = 20,000}

In our last example, we demonstrate time-accurate implicit solution of
transient flow around an SD7003 airfoil at Re = 20,000. The mesh is highly
resolved in the boundary layer, however, for this Reynolds number it is
coarser than the flow features in much of the domain and the computations
should therefore be considered an under-resolved ILES-type model
\cite{uranga11iles}.

\begin{figure}[t]
\begin{minipage}{.49\textwidth}
  \begin{center}
    \includegraphics[width=.9\textwidth]{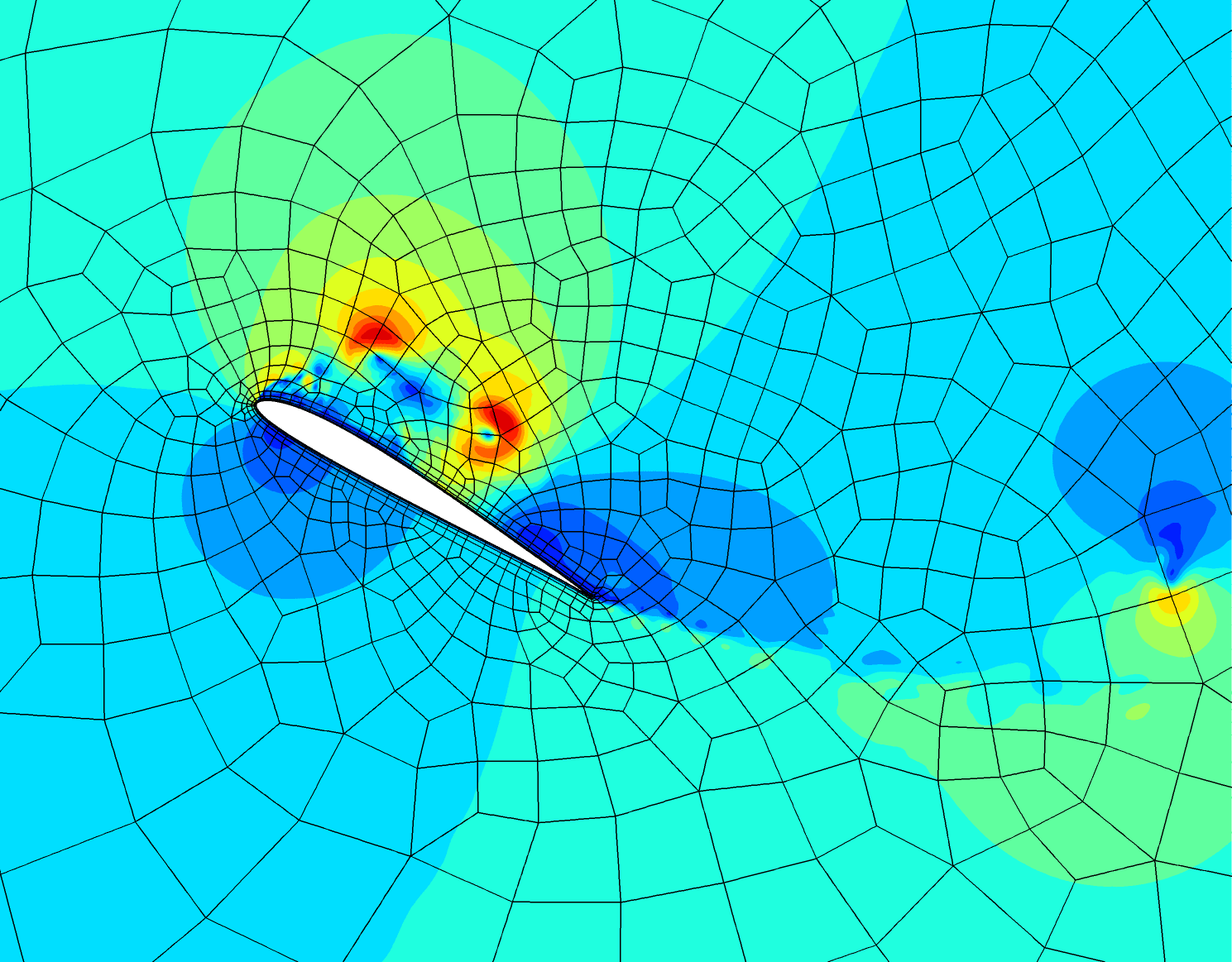} \\
  \end{center}
\end{minipage} \ %
\begin{minipage}{.49\textwidth}
  \begin{center}
    \includegraphics[width=.9\textwidth]{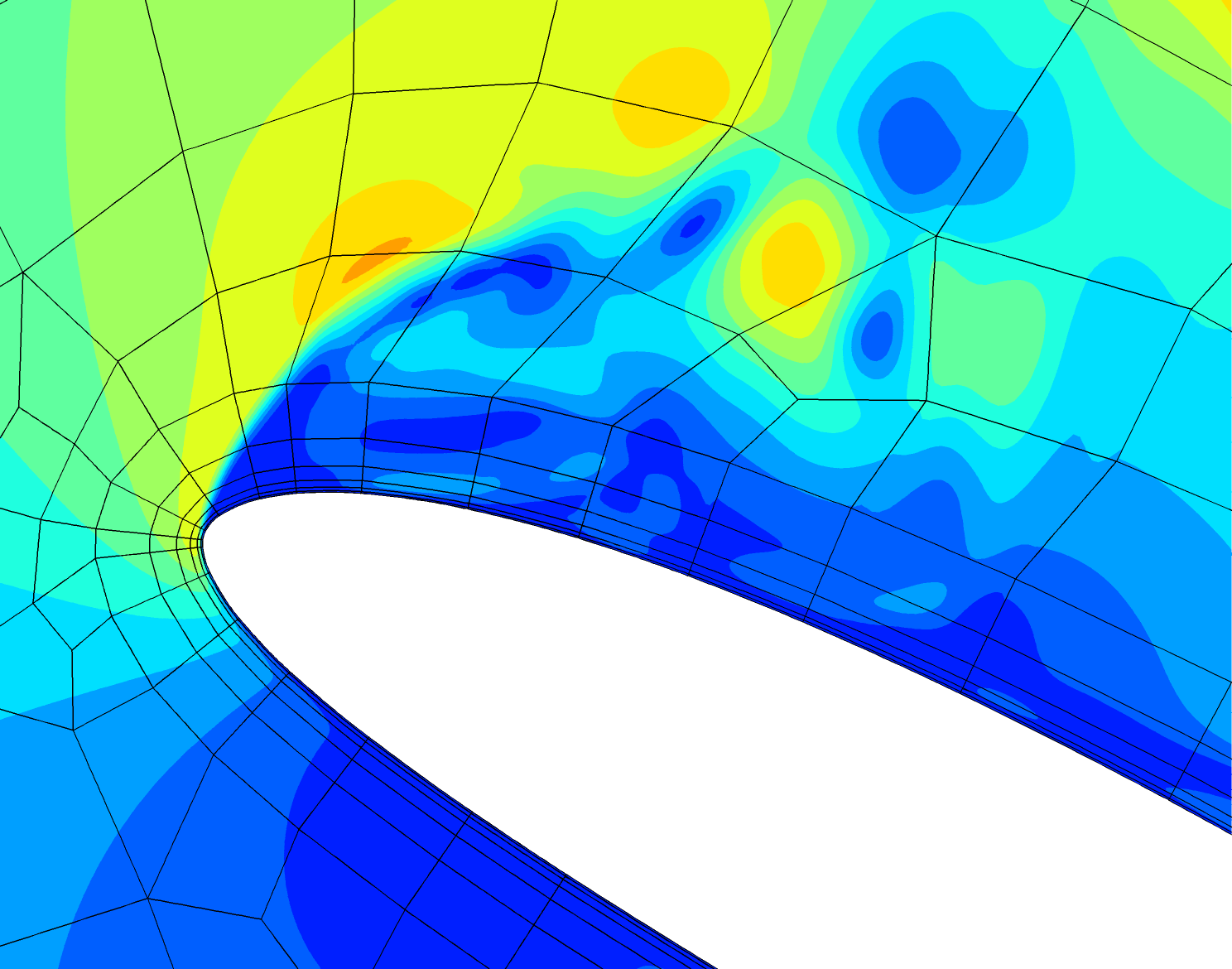} \\
  \end{center}
\end{minipage} \\
\begin{minipage}{.49\textwidth}
  \begin{center}
    \includegraphics[width=.95\textwidth]{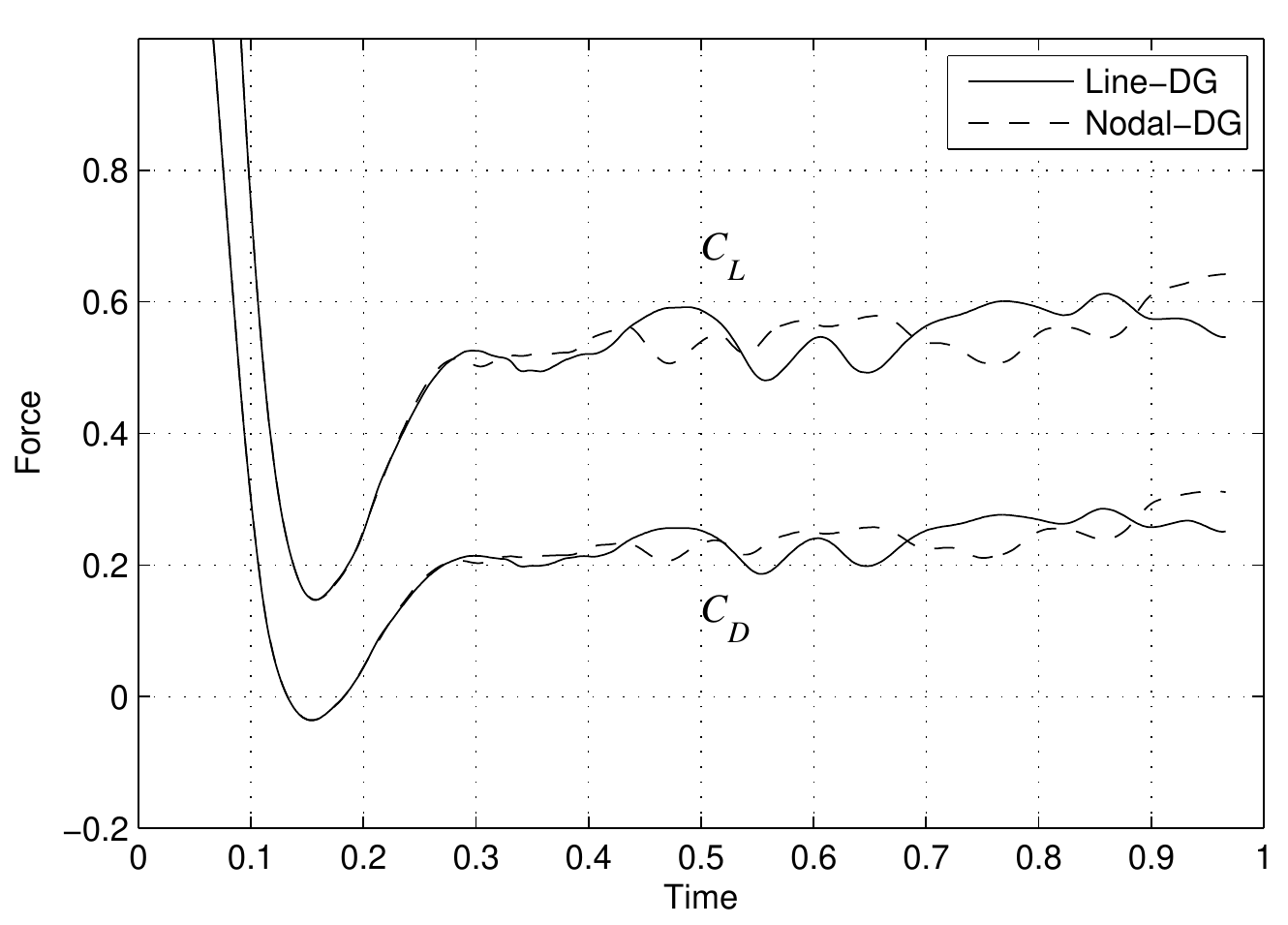} \\
  \end{center}
\end{minipage} \hfill
\begin{minipage}{.49\textwidth}
  \begin{center}
    \includegraphics[width=.95\textwidth]{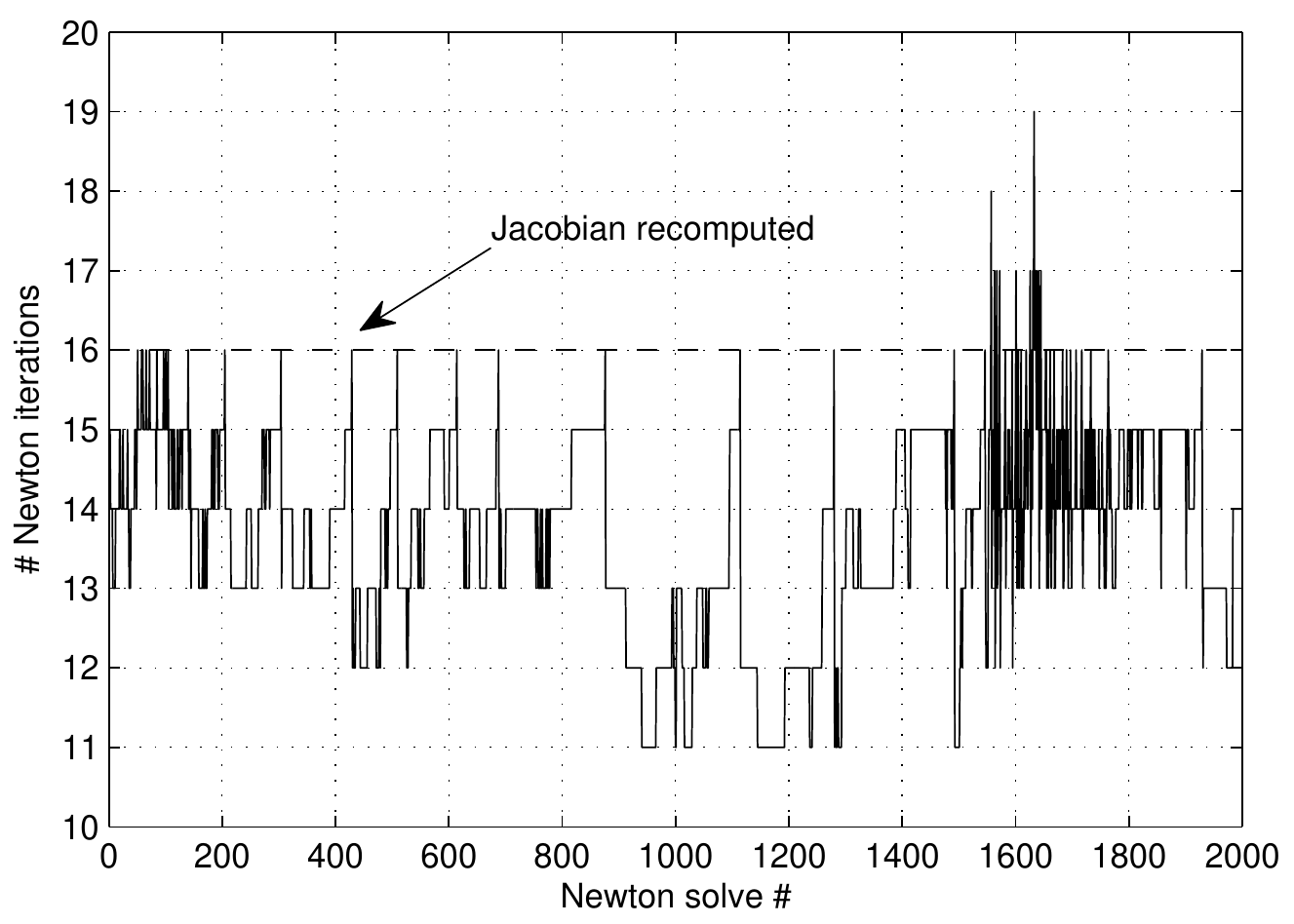} \\
  \end{center}
\end{minipage}
\caption{Implicit transient simulation of flow around an SD7003 airfoil, at 30
  degrees angle of attack, Reynolds number 20,000, and Mach number 0.1.
  Line-DG with polynomial degrees of $p=7$ is used for the spatial
  discretization, and a three stage, third order accurate DIRK scheme is used
  for time integration. The nonlinear systems are solved using a Newton-Krylov
  solver with re-used Jacobian matrices. The CFL number compared to the
  explicit RK4 method is about 250. The top figures show the computational
  mesh and color contours of the Mach number, from 0 to 0.3. The bottom left
  plot compares the lift and drag coefficient with a nodal DG scheme,
  indicating good agreement until the small differences between the schemes
  have grown enough to make the solutions hard to compare. The bottom right
  plot shows the performance of the nonlinear Newton solver, and in particular
  how it only recomputes the Jacobian matrices once in about every 100th
  solve.}
\label{sd1t}
\end{figure}

The Mach number is 0.1 and the angle of attack is 30 degrees to force flow
separation at the leading edge. We use the three-stage DIRK scheme
(\ref{dirk1}), (\ref{dirk2}), solved with Newton's method as described in
sections~\ref{sec:newtonkrylov} and~\ref{sec:reuse}. At each Newton step we
perform 5 GMRES iterations, and if the number of Newton iterations exceeds 15
we recompute the Jacobian matrices. Our computational mesh has 1122
quadrilateral elements with polynomial degrees $p=7$. The mesh and a solution
at the normalized time of $t=1.76$ are shown in the top plots of
figure~\ref{sd1t}.

We use the timestep $\Delta t=2\cdot 10^{-4}$, which is about 250 times larger
than the largest stable explicit RK4 timestep, yet small enough to accurately
capture most of the complex flow features. It is difficult to estimate the
accuracy in the simulation, due to the under-resolved nature of LES and the
high sensitivity of transitional flows. However, we have run the same problem
using a nodal DG code which has been tested against other simulations as well
as experiments \cite{uranga11iles}. The bottom left plot of figure~\ref{sd1t}
shows that the lift and drag forces on the airfoil agree well between the two
schemes, until small perturbations have grown enough to cause large differences
between the flows.

We also plot the performance of the Newton solver, in the bottom right of
figure~\ref{sd1t}. It shows how the number of Newton iterations per solve
remains fairly constant, and about once every 100th timestep it reaches 16
which forces a recomputation of the Jacobian matrix. The average number of
iterations per Newton solve for the entire simulation is 14. Because of the
sparsity of the matrices and the splitting (\ref{matvecsplit}), these
iterations are relatively inexpensive compared to residual evaluation.  Our
implementation is not optimized for performance, but the relative times for the
four operations (1) GMRES iteration, (2) residual evaluation, (3) Jacobian
evaluation, and (4) preconditioner factorization are roughly
1:4:40:16. Therefore, since we only perform 5 GMRES iterations per solve, we
spend about the same time in residual evaluation as in the linear solver. In
this sense, the solver is similar to an explicit scheme in that it spends a
large portion of its computational time in residual evaluations, which gives
benefits e.g. in the parallelization on new parallel multicore computer
architectures with limited memory bandwidth. The time for re-assembly and
factorization of the preconditioner is negligible, since they are only
performed once in about every 100th timestep, which corresponds to 1400
residual evaluations or 7000 GMRES iterations. We also point out that without
preconditioning about 20 times more GMRES iterations are required for the same
tolerance, showing that even our simple preconditioner makes a drastic
difference on the convergence.

\section{Conclusions}

We have presented a new line-based DG method for first and second-order
systems of equations. The scheme has a simple structure, with only
one-dimensional integrals and standard Riemann solvers applied point-wise.
Compared to the standard nodal DG method, this gives a simpler assembly
process and a fundamentally different sparsity structure, which we used to
develop efficient matrix-based implicit solvers. Compared to collocation based
methods such as the DG spectral element method, it uses fully consistent
integration along each coordinate-direction, and it slightly reduces the
connectivities to neighboring elements by the choice of solution nodes. We
showed that the accuracy of the discretizations are very similar to the
standard DG method, and we demonstrated a stiff LES-type flow simulation with
high-order DIRK time-integration with Newton-Krylov solvers and re-used
Jacobians.

A number of further developments are needed to make the scheme competitive for
real-world problems. We have not addressed the issue of nonlinear stability
for under-resolved features, including shock capturing, where approaches such
as artificial viscosity and limiting could be adopted. For the solvers, our
simple block-Jacobi preconditioner can be much improved upon, using e.g.
multigrid and ILU techniques. Finally, for large problems the implementation
needs to be parallelized, in particular for the new generation of multicore
and GPU chips where memory bandwidth is limited. Here the high sparsity of the
Line-DG scheme might have additional benefits over the standard nodal DG
scheme.

\section{Acknowledgments}

We would like to acknowledge all the valuable discussions about this work with
Jaime Peraire, Luming Wang, and Bradley Froehle, the suggestions from the
reviewers, as well as the generous support from the AFOSR Computational
Mathematics program under grant FA9550-10-1-0229, the Alfred P. Sloan
foundation, and the Director, Office of Science, Computational and Technology
Research, U.S. Department of Energy under Contract No. DE-AC02-05CH11231.

\bibliographystyle{model1-num-names}
\bibliography{linedg}

\end{document}